\RequirePackage{fix-cm}

\documentclass[smallextended]{svjour3bis}       
\smartqed  
\usepackage{graphicx}
\usepackage{amsmath, amssymb, amsfonts}
\usepackage{booktabs}
\usepackage[authoryear]{natbib}

\usepackage{dsfont}
\usepackage{relsize}
\usepackage{array,multirow,makecell} 
\usepackage{hhline}
\usepackage{bm}
\usepackage{placeins}
\usepackage{rotating}
\usepackage{graphicx}
\RequirePackage[colorlinks,citecolor=blue,urlcolor=blue]{hyperref}
\DeclareMathOperator*{\argmin}{argmin}%

 \journalname{Submitted to Scandinavian Journal of Statistics}

\def \V{\mathbb V}
\def \R {\mathbb R}

\def \beg{\begin{eqnarray}}
\def \en{\end{eqnarray}}
\def \be*{\begin{eqnarray*}}
	\def\e*{\end{eqnarray*}}
\def \beq*{\begin{equation*}}
\def\eq*{\end{equation*}}
\def \di{\displaystyle}
\def\bit{\begin{itemize}}
	\def \eit{\end{itemize}}
\def \E{\mathbb E}
\def \P{\mathbb P}
\def \N{\mathbb N}

\def \w{\widehat}
\def \t{\tilde}

\def \m {\mathbf}

\graphicspath{{./Img/}}

\begin{document}

\title{Nonparametric estimation of copulas and copula densities by orthogonal projections
}

\titlerunning{Nonparametric estimation of copulas}        

\author{{\sc   Yves I. Ngounou Bakam$^{1}$ }        \and
        {\sc Denys Pommeret$^{1,2}$ }     
}


\institute{{\bf Correspondence}
\at Denys Pommeret, ISFA, 50, Avenue Tony Garnier
69366 LYON CEDEX 07, France.
 \\  \email{denys.pommeret@univ-amu.fr}
\and
$~^{1}$  CNRS, Centrale Marseille, I2M, Aix-Marseille Univ., Marseille, France 
           \\
         $~^{2}$  Univ Lyon, UCBL, ISFA LSAF EA2429, F-69007, Lyon, France
}

\date{}

\maketitle

\vspace*{-2cm}
\begin{abstract}
In this paper we study nonparametric estimators of copulas and copula densities.
We first focus our study on a copula density estimator based on polynomial  orthogonal projections  of
the joint density. 
A new copula estimator is then deduced simply by integration. 
Both estimators are based on a characteristic sequence of the copulas, which we refer to as the {\it copula coefficients}. 
Their asymptotic properties are reviewed
and a data driven selection of the number of coefficients allows these estimators to be well adapted to all types of copulas. In particular, this selection procedure perfectly  detects the cases of independence.
An intensive  simulation study
 shows the very good performance of both copulas and  copula densities estimators in comparison to a large panel of competitors. 
Applications in actuarial science illustrate this approach.
\keywords{Copula estimators \and Copula coefficients \and 
Minimax theory  \and Orthogonal Legendre polynomials}
\end{abstract}

\section{Introduction}

Consider a $d$-random vector $\mathbf{X} = (X_1, \cdots , X_d)^T$ with joint distribution function $H$ and marginal distribution functions $F_1,\cdots , F_d$, that we assumed to be continuous.
 According to Sklar's Theorem \citep{Sklar}, there exists a unique $d$-variate
 function $C$ such that
 \beg
 H(x_1,\cdots, x_d) &=& C(F_1(x_1),\cdots , F_d(x_d)).
 \en
  The function $C$ is called the copula  associated to $\mathbf{X}$. The copula is a joint  distribution
 function on $I^d=[0, 1 ]^d$, with uniform margins and satisfying
$C(u_1,\cdots, u_d) = H(F_1^{-1} (u_1), \cdots, F_d^{-1} (u_d))$, where, for $j = 1,\cdots, d$, $F^{-1}_j(u_j) = \inf\{x_j;  F_j(x_j) \geq  u_j\}$,
is the quantile function of $F_j$.
Assuming that for $j = 1,\cdots, d$, $F_{j}$ is differentiable, we can express the joint density $h$ of $\mathbf{X}$ (with respect to the Lebesgue measure on $I^{ d}$) as
 \be*
\label{cop1}
h(x_1,\cdots,x_d)  = c\left(F_{1}(x_1),\cdots,F_{d}(x_d)\right)\prod_{j=1}^{d}f_{j}(x_j),
\e*
where for $j=1,\cdots,d$, $f_{j}$  is the marginal density of ${X_j}$ and where
\be*
\label{density}
c & =  & \frac{\partial^{d} C}{\partial x_1\cdots \partial x_d },
\e*
is called the  copula density of $\mathbf{X}$.

Copulas and copulas density have a large spectra of applications as described for instance in  \citet{Joe}.
They are largely used
 as a tool to identify a wide variety of properties such as tail-dependencies, heavy-tail, asymmetries or the heavy-tail behaviour.
 Copula estimation has given rise to numerous works.
 In parametric approaches
 maximum likelihood estimation is commonly used, with a general two-stage procedure as  in  \citet{Ko}, or penalized likelihood as in \citet{Qu2012}, and more recently for Archimax copulas in \citet{Chatelain} . 
 Such parametric methods can
suffer from restrictive specification errors with few  parameters.  Nonparametric copula estimations generally offer a more
flexible alternative 
 and in such framework copulas density estimators have been studied with different methods: empirical processes \citep[see][]{deheuvels1979fonction},  kernel estimators  \citep[see for instance][]{Geenens, omelka2009improved} where the authors proposed a robust method,
B-spline estimators \citep{Kauermann}, 
Bernstein polynomials  \citep{bouezmarni2010asymptotic, bouezmarni2013bernstein}, 
linear wavelet estimators  \citep[see][]{Genest}, nonlinear wavelet estimators \citep[see][]{Autin} or Legendre multiwavelet estimators \citep[see][]{chatrabgoun2017legendre}.

In this paper we first propose to revisit the nonparametric method
for estimating copula densities by considering an
orthogonal shifted Legendre polynomials expansion.
In this sense, our approach lies between the Legendre
multiwavelet procedure  and
the Bernstein method.
More precisely, the basic idea developed here is to consider the transformations $U_j=F_j(X_j)$, for $j=1,\cdots, d$, which yield  to a vector of uniform random variables  denoted by $$\mathbf{U}=(U_1, \cdots, U_d)^T.$$ Its   joint distribution function has the form
\begin{align*}
    H_{\mathbf U}(u_1,\cdots, u_d)
    & =
    H(F_1^{-1}(u_1), \cdots, F_d^{-1}(u_d)),
\end{align*}
and clearly, $\mathbf U$ and $\mathbf X$ have the same structure of dependence with the same copula. We deduce that
\begin{align}
\label{densite2}
h_{\mathbf U}(u_1,\cdots, u_d) &
=
c(u_1,\cdots, u_d),
\end{align}
where $h_{\mathbf U}$ denotes the joint density of $\mathbf U$ with respect to the Lebesgue measure on $[0,1]^{ d}$.
This basic equality is the start of the construction of the copula density estimator, expressing $c$ in a basis of orthogonal polynomials.   In addition,
a new copula estimator is derived by simply integrating the polynomial expansion. Both estimators are based on a characteristic sequence of the copulas, which we refer to as the {\it copula coefficients}. 
The properties of these estimators are studied: the copula density estimator satisfies the uniform margins property (see Proposition
\ref{pro_uniform}) that almost all
nonparametric estimators in the literature suffer.
Various asymptotic properties are reviewed for both estimators.
 We also provide a  functional class for which the construction of the copula estimator is optimal in the minimax  
 sense. We propose a data driven method of selection  for the number of projections. In the case of independence this criteria can find the exact form of the copula and  our estimator then seems to greatly surpass all the other competitors. 
 Moreover, a relation between the Spearman's rho $\rho_C$ and the copula coefficients   is highlighted, yielding to a new estimator of $\rho_C$.
 Finally, both  estimators seem to outperform  all others  known in the statistical literature for a wide
spectrum of scenarios.
Moreover, these estimators  are very simple and easy to
implement and their execution time is very fast.
Their interest is also demonstrated on a real dataset.

The paper is organised as follows. In section 2 we develop the methodology of the estimation procedure. Section 3 is devoted to the elementary and asymptotic properties with minimax 
results. A selection method of
the number of expansion components and a numerical comparison with other nonparametric
estimators are proposed in Section 4.  Section 5 is devoted
to two  applications on real actuarial dataset. Finally in Section
6 we discuss further leads of research connected with copula estimation and we explain why our approach is efficient even without square integrability assumption. All proofs are relegated to the Appendix.

\vspace{0.5cm}

\section{Construction of the estimators}
\label{methodo}
Let us denote by $\mu$ the uniform measure on $\di I=[0,1]$  and by  
$\{Q_{m}; m\in\N\}$ an orthonormal basis
of shifted Legendre polynomials satisfying
\be*
\di\int_{I} Q_{m}(x) Q_{k}(x) \mu(dx) &=& \delta_{mk},
\e*
with $\delta_{mk}=1$ if $m=k$ and $0$ otherwise.
The orthonormal shifted Legendre polynomials $\di Q_{m}$ are defined on $\di I$ by
\be*
\di Q_{m}(x)=\sqrt{2m+1}L_{m}(2x-1),
\e*
where $L_{m}$ are the well known classical Legendre polynomials defined by
\be*
L_0=1;  L_1(x) = x; &{\rm and } &
(m+1)L_{m+1}(x) = (2m+1)xL_m(x)-mL_{m-1}(x).
\e*
Characterizations and properties of Legendre polynomials can
be found in \citet{abramowitz1970handbook}. 
For all $\mathbf{x}=(x_1,\cdots,x_d)^T \in I^d$ and for all $\mathbf{m}=(m_1,\cdots, m_d)^T \in \N^d$, we define
\be*
\mathbf{Q}_{\mathbf m}(\mathbf{x}) & = & \prod_{j=1}^{d}Q_{m_j}(x_j),
\e*
and we write
\beg
\label{coef}
\rho_{\mathbf m} & = &
\E\left(\prod_{j=1}^{d}Q_{m_j}(F_{j}(X_j))\right).
\en
The following assumption will be used  throughout the paper  
and says 
that the copula density $c$ 
belongs to $\mathbb L^2 ([0,1]^d)$
\beg
\label{condition}
\di \lVert c \rVert_2^2 = \sum_{{\mathbf m}\in\N^d}\rho_{\mathbf{m}}^2 <\infty,
\en
where $\lVert . \rVert_2$ denotes the norm with respect to the Lebesgue measure on $[0,1]^{ d}$. 

Assumption (\ref{condition}) 
is obviously satisfied for any bounded copula density,
for instance the Farlie-Gumbel-Morgenstern and Frank copula densities (see \cite{Nelsen} for
definitions). \cite{beare2010copulas}  showed that $\di \lVert c_\rho \rVert_2^2= (1-\rho^2)^{-1/2}$ for the standard bivariate Gaussian copula density with correlation coefficient $\rho \in (0;1)$.  
 Moreover, \cite{beare2010copulas} noted that in the bivariate case, copulas associated to Lancaster type distributions \citep{lancaster} satisfied (\ref{condition}). This is the case for bivariate gamma, Poisson,
binomial and hypergeometric distributions, and for the compound correlated bivariate Poisson distribution (see for instance \cite{hamdan}).  
However, copulas exhibiting lower or upper tail dependence (in the sense of \cite{McN05}) do not have square integrable density.  In particular, the Gumbel, Clayton, and t-copulas all have upper or lower
tail dependence and then do not satisfy condition (\ref{condition}). 
 In Appendix \ref{appendL2}  we  discuss the possibility to modify such  copulas with a shrinkage factor to be square integrable and we illustrate such a  procedure.   
 This explains why our method can work so well even if (\ref{condition}) is not satisfied. But it is important to note that our framework is nonparametric and that the possible family of the copula is  unknown.

\begin{proposition}
\label{prop0}
Assume that 
(\ref{condition}) holds. 
Then we have
\beg
\label{classic}
 c(u_1,\cdots,u_d)  &=&
\di\sum_{\mathbf{m}\in\N^d}\rho_{\mathbf m}\mathbf{Q}_{\mathbf m}({\mathbf u})
\\
\label{cop2}
C(u_1,\cdots,u_d) & = & \di \sum_{{\mathbf m}\in\N^d}
\rho_{{\mathbf m}}\prod_{j=1}^{d}\di\int_{0}^{u_j}Q_{m_j}(x_j) \mu(dx_j).
\en
\end{proposition}

It is important to note that the sequence $(\rho_{{\mathbf m}})_{{\mathbf m
}\in \N^d}$  characterizes the copula. In this way it will be referred to as the \emph{copula coefficients}.
Since $Q_0=1$ we have $\rho_{\mathbf 0}=1$. The particular case $\rho_{\mathbf m}=0$ for all $\mathbf m\neq 0$ coincides with the independent case. As seen in (\ref{coef}), the sequence $(\rho_{{\mathbf m}})_{{\mathbf m
}\in \N^d}$ contains all the polynomial correlations between the marginal uniform random variables.
In the bivariate case, that is when $d=2$,  the  element $\rho_{\mathbf m}$ simply expresses the correlation between $Q_{m_1}(F_{1}(X_1))$ and $Q_{m_2}(F_{2}(X_2))$.
In the general case, by orthogonality, we have $\E\left(Q_{m_j}(F_{j}(X_j))\right)=0$ for all $m_j>0$ and then  $\rho_{\mathbf m}=0$ as soon as only one component of $\m m$ is zero.
Moreover, if for some integer $j\in\{1,\cdots,d\}$,   $X_j$ is independent to all over variables $X_i$, $i\neq j$, then $\rho_{\mathbf m}=0$
as soon as $m_j>0$.
Then the copula coefficients can be used as an indicator of   independence between the components of $\mathbf X$. In this sense, we also exhibit a link with the Spearman's rho in Section \ref{some properties}.

For any  positive  integer vector ${\mathbf N}=(N_1,\cdots,N_d)^T$ we  define the following ${\mathbf N}$-th order approximations: 
\be*
c^{[\mathbf{N}]}(\mathbf{u}) & = &
\sum_{{\mathbf m}\leq \mathbf{N}}\rho_{{\mathbf m}}\prod_{j=1}^{d}Q_{m_j}(u_j)
\\
C^{[\mathbf{N}]}(\mathbf{u}) & = &
\di \sum_{{\mathbf m}\leq \mathbf{N}}
\rho_{{\mathbf m}}\prod_{j=1}^{d}\di\int_{0}^{u_j}Q_{m_j}(x_j)\mu(dx_j),
\e*
where the inequality ${\mathbf m}\leq {\mathbf N}$ means that $m_j\leq N_j$ for all $j=1,\cdots, d$. We write $\mathbf m \not\leq \mathbf N$ when $\mathbf m$ does not satisfy this inequality.
If we  observe a $n$-sample
$\mathbf{X}_{1},\cdots, \mathbf{X}_{n}$,  of iid random data,  with $\mathbf{X}_i=(X_{i1},\cdots, X_{id})^T$ having joint distribution function $H$,  
 then we can estimate the quantity $\rho_{\mathbf m}$ by
\be*
\label{estimrho}
\w{\rho}_{{\mathbf m}} & = &
\left\{
\begin{array}{cc}
\di 1     & {\rm if \ }\mathbf m=\mathbf 0,  \\
\di 0     &  {\rm if  \ exactly \ }d-1 {\rm \ components \ of \  } \mathbf m {\rm \  are  \ zero},
\\
\di \frac{1}{n}\di\sum_{i=1}^{n}
\prod_{j=1}^{d}Q_{m_j}(\w F_{j}(X_{ij})),
& {\rm else},
\end{array}
\right.
\e*
where $ 
\di\w{{F}}_j({x}) = \frac{1}{n}\sum_{i=1}^{n}\mathds{1}(X_{ij}\leq x)$.

\begin{remark} 
Since the marginal distributions are continuous, the ties occur with probability zero. So we can also simply write
\be*
\w{\rho}_{{\mathbf m}} & = & \frac{1}{n}\di\sum_{i=1}^{n}\prod_{j=1}^{d}Q_{m_j}\Big(\frac{R_{ij}}{n}\Big),
\e*
where 
$R_{ij}$ denotes the rank of $X_{ij}$.
\end{remark}

A $\mathbf N$-th order nonparametric estimator of the copula density $c$ is given by
\beg
\label{copknown}
\di \w{c}^{[{\mathbf N}]}(u_1,\cdots,u_d) &=& \sum_{{\mathbf m}\leq {\mathbf N}}\w{\rho}_{{\mathbf m}}\prod_{j=1}^{d}Q_{m_j}(u_j).
\en
By integration, we get a very simple $\mathbf N$-th order  nonparametric  estimator of the copula function  as follows
\beg
\label{estcop3}
\di \w{C}^{[{\mathbf N}]}(u_1,\cdots,u_d) = \sum_{{\mathbf m}\leq {\mathbf N}}\w{\rho}_{{\mathbf m}}\prod_{j=1}^{d}\displaystyle {\int_{0}^{u_j}Q_{m_j}(x_j) \mu(dx_j)}.
\en
\begin{remark}
\label{remark2}
In the particular case where  the margins $F_1,\cdots, F_d$ are known, the  pseudo-estimator of the copula  coefficients is given by
\be*
\label{estimrho2}
\tilde{\rho}_{{\mathbf m}} & = &  \frac{1}{n}\di\sum_{i=1}^{n}
\prod_{j=1}^{d}Q_{m_j}(F_{j}(X_{ij})),
\e*
and the associated copula and copula density pseudo-estimators are
\be*
\label{copknown2}
\di \tilde{c}^{[{\mathbf N}]}(u_1,\cdots,u_d) &=& \sum_{{\mathbf m}\leq {\mathbf N}}\tilde{\rho}_{{\mathbf m}}\prod_{j=1}^{d}Q_{m_j}(u_j),
\e*
\end{remark}
\beg
\label{estcop32}
\di \tilde{C}^{[{\mathbf N}]}(u_1,\cdots,u_d) = \sum_{{\mathbf m}\leq {\mathbf N}}\tilde{\rho}_{{\mathbf m}}\prod_{j=1}^{d}\displaystyle{\int_{0}^{u_j}Q_{m_j}(x_j) \mu(dx_j)}.
\en
\begin{remark} 
 We can see that  
$\w c^{[\m 0]} =1$ which coincides exactly  with the
copula density  in the independent case.  
In such case, we will observe this phenomenon in our simulation study where the copula and its estimator are very often exactly the same. 
\end{remark}

\section{Some properties of the  estimators}\label{some properties}
Write $\bm \mu^d = \mu \times \cdots \times \mu$ the product of the uniforme measure $d$ times. 
\subsection{Elementary properties}
\begin{proposition}\label{bona_fide}
Let $\mathbf{u} \in I^d$ and fix $\m N \in \N^d$. The copula estimator $\di \w{C}^{[{\mathbf N}]}$ given by  (\ref{estcop3}) satisfies the following properties
\begin{itemize}
\item[i)] If at least one coordinate of $\mathbf{u} $ is zero then $\di \w{C}^{[{\mathbf N}]}({\mathbf u})= 0$.
\item[ii)] If  $\di \mathbf{u}=(1, \cdots,1,u_i,1, \cdots,1)$, then $\w{C}^{[{\mathbf N}]}({\mathbf u})=u_i$.
\item[iii)]$\di \w{C}^{[{\mathbf 0}]}({\mathbf u})=\prod_{j=1}^{d}u_j$. 
\end{itemize}
\end{proposition}

\begin{proposition}\label{pro_uniform}
The copula density estimator $\di \w{c}^{[{\mathbf N}]}$ given by  (\ref{copknown}) satisfies the following properties
\begin{itemize}
\item[i)] 
$\di\int_{I^{d}}\w{c}^{[{\mathbf N}]}({\mathbf u})\bm{\mu}^d (d{\mathbf u})=1$.
\item[ii)]
For all $u_j \in I$, $j=1,\cdots,d$, writing  \,$\di {\mathbf x_{-j}}=(x_1,\cdots,x_{j-1},u_j,x_{j+1},\cdots,x_d)$, we have
\be*
\di \int_{I^{d-1}}\w{c}^{[{\mathbf N}]}({\mathbf x_{-j}})\bm{\mu}^{d-1}(\mathrm{d}{x_1 \cdots \mathrm{d} x_{j-1} \mathrm{d} x_{j+1} \cdots \mathrm{d} x_d}) = 1 .
\e*
\item[iii)] $\di \w{c}^{[{\mathbf 0}]}({\mathbf u})=1$, for all ${\mathbf u} \in I^d$.
\end{itemize}
\end{proposition}
Let us recall that for any continuous bivariate random variable $(X_1,X_2)$ with copula $C$,  the Spearman's rho  can be express as 
(see \cite{Nelsen}):
\be*
\di \rho_C = 12\int_{0}^{1}\int_{0}^{1}C(u,v)dudv-3.
\e*
The result below gives a  simplified new  expression of $\di \rho_C$ in terms of copula coefficients.  

\begin{proposition}\label{theo spearman}
Let $(X_1,X_2)$ be a  continuous bivariate random variable with copula $C$. Then the Spearman's rho coincides with the first copula coefficient defined by (\ref{coef}), that is:   
\be*
  \di \rho_C =\rho_{11}.
\e*
\end{proposition}
We can immediately deduce an  estimator of the Spearman's rho as follows:
\be*
 \di \hat{\rho}_C =
 \di \hat{\rho}_{11}
 = \frac{3}{n}\sum_{i=1}^{n}\left(2\w U_{i1}-1\right)\left(2\w U_{i2}-1\right).
\e*
This estimator is new and could  be compared to the ones given in \cite{perez2016note} but such a study exceeds the scope of this paper.

\subsection{Asymptotic properties}
\begin{proposition}\label{prop 1}
Fix $\mathbf N \in \N^d$,  independent of $n$. For all $\mathbf{u}\in I^d$, we have
\begin{itemize}
\item[i)]
\be*
\w{c}^{[{\mathbf N}]}({\mathbf u}) &=&\tilde{c}^{[{\mathbf N}]}({\mathbf u})+\mathcal{O}_{\P}\big(\sqrt{n^{-1}}\big).
\e*

\item[ii)]
\be*
\w{C}^{[{\mathbf N}]}({\mathbf u})&=&\tilde{C}^{[{\mathbf N}]}({\mathbf u})+ \mathcal{O}_{\P}\big(\sqrt{n^{-1}}\big)
\e*

\end{itemize}
\end{proposition}
 For any integer $\mathbf{N} \in \N^d$, the notation $\mathbf{N} \to \infty$ means $N_j \to \infty$ for all $j=1,\cdots, d$ and the $\di \max_{j=1,\cdots,d}(N_j)$ will be denoted by $\di N_{max}$. 

We now consider the Mean Integrated Squared Error (MISE) as a rule criteria to decide which degree of approximation we use.  We write
\be*
MISE(\w{c}^{[\mathbf{N}]})&=&\di\E \| \w{c}^{[\mathbf{N}]}-c\|_{2}^2=\E \int_{I^d}\left(\w{c}^{[{\mathbf{N}}]}(\mathbf{u})-c(\mathbf{u})\right)^2 \bm \mu^d(d\mathbf{u}).
\e*
\begin{proposition}\label{prop 4}
 Assume that (\ref{condition}) holds. For all ${\mathbf N} \in\N^d$ we have
 \be*
\label{mise}
MISE(\w{c}^{[\mathbf{N}]})
 &=& \E\left(\sum_{\mathbf{m} \leq \mathbf{N}}(\w{\rho}_{\mathbf{m}}-\rho_{\mathbf{m}})^2\right) +  \di \sum_{\mathbf{m} \not\leq \mathbf{N}}\rho_{\mathbf{m}}^2. 
\e*
\end{proposition}
Thus $\mathbf{N}$ represents a smoothing parameter which controls the trade-off between the bias-squared and the variance. \begin{corollary}\label{relation_mise}
Fix $N\in\N^d$, independent of $n$. Then we have 
\be*
\di MISE(\w{c}^{[\mathbf{N}]})=MISE(\widetilde{c}^{[{\mathbf{N}}]})+\mathcal{O}(n^{-1})
\e*
\end{corollary}

From now on, assume that $\di {\mathbf{N}}={\mathbf{N}}(n)$ and consider the assumption
\begin{align*}
{(\mathbf H)}  & \hspace{3cm}
 N_{max}^{2d+4} =
o\left( n \right).\hspace{5cm}
\end{align*}

\begin{corollary} \label{corol 1}
Assume that (\ref{condition}) and {\bf (H)} hold. 
The copula density estimator $\w{c}^{[{\mathbf{N}}]}$ is asymptotically consistent in the integrated mean squared sense, that is:
\be*
MISE(\w{c}^{[\mathbf{N}]})
\xrightarrow[\enskip n\longrightarrow \infty\enskip]{} 0.
\e*
\end{corollary}

\begin{proposition}\label{prop 2}
Assume that (\ref{condition}) and {\bf (H)} hold. For any $\mathbf{u} \in I^d$, 
we have as $n \to \infty$
\be*
\di \E\left(\w{c}^{[{\mathbf N}]}(\mathbf{u})\right) \longrightarrow  c(\mathbf{u}) 
{\rm \ and \ }
\E\left(\w{C}^{[{\mathbf N}]}(\mathbf{u})\right) \longrightarrow  C(\mathbf{u})
\e*
\end{proposition}
\subsection{Minimax results}

We recall here the minimax procedure in the spirit of \citet{Lepski} or  \citet{Cohen}.  
of a given  estimator $\w{c}$
over this set 
by its maximum risks (or worst case risk) defined as follows
\be*
\di \mathcal{R}( \mathcal{F})=\sup_{c\in \mathcal{F}}\E\|\w{c}-c\|^2_2,
\e*
that we compare to a benchmark which is the following minimax risk
\be*
 \di \mathcal{R}^{\star}( \mathcal{F}, {\cal C})=\inf_{\w{c}\in {\cal C}}\sup_{c\in \mathcal{F}}\E\|\w{c}-c\|^2_2,
\e*
 where $\cal C$ denotes the  class of the  estimators $\w c$. In our case we consider the class of estimators satysfying   (\ref{copknown})  which depends of the approximation degree $\m N$. Then we simply write 
 \be*
 \di \mathcal{R}^{\star}( \mathcal{F})=\inf_{{\m N}}\sup_{c\in \mathcal{F}}\E\|\w{c}^{[\m N]}-c\|^2_2. 
\e* 

 Denoting by $n$ the sample size, a sequence $(v_n)$ is said to be an optimal rate convergence in the minimax sense for $\mathcal{F}$ if there exits two constants $M_1$ and $M_2$ such that
\be*
\di M_1 v_n\leq \mathcal{R}^{\star}( \mathcal{F})\leq M_2 v_n,
\e*
 and the estimator $\w{c}$  is said minimax optimal if there exists a constant $ \Upsilon$ such that 
\be*
\mathcal{R}( \mathcal{F})\leq \Upsilon v_n.
\e*
 For more details about minimax theory we alos refer  to \citet{Tsybakov} and 
 \citet{efromovich2008nonparametric}.

Under (\ref{condition}) the copula density $c$ is characterized by its sequence of copula coefficients  $(\rho_{{\mathbf{m}}})_{{\mathbf{m}}\in \N^d}$, and  for  $\boldsymbol{\beta}=(\beta_1,\cdots,\beta_d)\in (\R_{+}^*)^{d}$,  and $0<L<\infty$, we denote by $\mathlarger{\mathcal{F}}_{\boldsymbol{\beta}}(L)$ the ellipsoid space  defined by:
\be*
\mathlarger{\mathcal{F}}_{\boldsymbol{\beta}}(L)=\left\{ c \in \mathbb{L}^{2}([0,1]^d)
~/\di \sum_{\mathbf{m} \in \N^d}{\rho}_{\mathbf{m}}^{2}\left(1+\sum_{i=1}^{d}m_{i}^{\beta_i}\right) < 
L 
\right\}.
\e*

As discussed in Appendix \ref{appendL2} there are some copulas families which do not satisfy assumption (\ref{condition}) and then do not belong to  $\mathlarger{\mathcal{F}}_{\boldsymbol{\beta}}(L)$. This is the case for Clayton, Gumbel or Student copulas. 
However, such copulas can be approximated as closely as we want by a function of this space as  illustrated in Appendix \ref{appendL2}.

The fonctionnal space $\mathlarger{\mathcal{F}}_{\boldsymbol{\beta}}(L)$ with assumption $c \in \mathlarger{\mathcal{F}}_{\boldsymbol{\beta}}(L)$ allows to control the bias term of the risk. We observe that, on the balls  $\mathlarger{\mathcal{F}}_{\boldsymbol{\beta}}(L)$, we have 
\beg\label{bias_term}
\|c^{[\m N]}-c\|^2_2 \leq L\sum_{j=1}^{d}N_{j}^{-\beta_j}. 
\en
\begin{remark}
\label{rem3}
Usually multivariate regularity spaces are defined
by 
\be*
\mathlarger{\mathcal{\t{F}}}_{\boldsymbol{\beta}}=\left\{ c \in \mathbb{L}^{2}(I^d)
~/\di \sum_{\mathbf{m} \in \N^d}{\rho}_{\mathbf{m}}^{2}\prod_{i=1}^{d}m_{i}^{\beta_i} < \infty
\right\},
\e*
 and equation (\ref{bias_term}) is also satisfied on the balls
\be*
\mathlarger{\mathcal{\t{F}}}_{\boldsymbol{\beta}}(L)=\left\{ c \in \mathbb{L}^{2}(I^d)
~/\di \sum_{\mathbf{m} \in \N^d}{\rho}_{\mathbf{m}}^{2}\prod_{i=1}^{d}m_{i}^{\beta_i} < L
\right\}.
\e* 
Note that all results below can also established on such  spaces $\mathlarger{\mathcal{\t{F}}}_{\boldsymbol{\beta}}(L)$. 
\end{remark}



\begin{proposition}\label{theor1} 
Assume that  (\ref{condition}) holds.  Let $\di b^{-1}=\sum_{i=1}^{d}\frac{1}{\beta_i}$ and assume that the integer vector $\m N$  satisfies 
\begin{align*}
{(\mathbf H')}  & \hspace{3cm}
N_{max}=\lfloor n^{\frac{1}{2d+b+4}}\rfloor,
\hspace{5cm}
\end{align*} 
where 
$\lfloor.\rfloor$  denotes the integer  part.
Then for all sample size $n$ we have 
\be*
\sup_{\di c\in {\mathcal{F}}_{\boldsymbol{\beta}}(L)}\E
\| \w{c}^{[{\mathbf{N}}]}-c\|^{2}_{2}& < & \alpha n^{-\frac{b}{b+2d+4}},
\e*
where $\alpha = 8(d)^2 ({3})^{d} + {2({3})^{d}}+Ld$. 
\end{proposition}
 Clearly, the previous condition {\bf (H')} is stronger than {\bf (H)}. 
 We known that in the   independent case (see Remark \ref{rem3}), if $\m N=\m 0$ we have  
 $\w c^{[\m 0]} = c$. We now exclude the case $\m N =\m 0$ and we 
 propose an upper bound of the minimax risk when $ \di \mathcal{R}^{\star}( \mathcal{F})=\inf_{\m N \neq \m 0}\sup_{c\in \mathcal{F}}\E\|\w{c}-c\|^2_2$.  
%
%
%
%
\begin{proposition} \label{coellip} 
Assume that (\ref{condition}) holds and
$\m N \neq \m 0$. Let  $v_n=n^{-\frac{1}{2d+b+4}}$. Then 
\be*
 \di \mathcal{R}^{\star}\left( \mathcal{F}_{\boldsymbol{\beta}}(L)\right)\leq K v_n, 
\e*
where 
$K=\frac{Ld(2d+b+4)}{2d+4}\left(\frac{Lbd}{\eta(2d+4)}\right)^\frac{-b}{2d+b+4}$, with $\eta = 8(d)^2 ({3})^{d} + {2({3})^{d}}.$


\end{proposition}


%


\section{Numerical analyses and comparisons}\label{simul-studies}
In this section, we present
simulation results which demonstrate the performance of our approach compared with a number of recent alternative estimators.
All computations were performed using the R software.

Since the results depend of the degree of the approximations we first present a data driven  method to select $\m N$.

\subsection{Data driven degree  selection}
\label{subsection 3.3}
We propose to use a data-driven procedure based on the Least-Squares Cross-Validation (LSCV) to select the optimal parameter $\di \w{\mathbf{N}}_{opt}$.
The LSCV procedure has been introduced by 
\citet{rudemo} and \citet{bowman} to
select the smoothing bandwidth for Kernel density estimation
and it has been  adapted to orthogonal series estimators by
 \citet{taylor}.
In the general case, the smoothing parameter is the minimizer of the following function
\be*
\di LSCV({\mathbf N})= \int_{I^{d}} \left(\w{c}^{[{\mathbf N}]}({\mathbf u})\right)^2 d{\mathbf u} -\di \frac{2}{m}\di \sum_{i=1}^{m}\w{c}^{[{\mathbf N}]}_{-i}\left(F_{1}(X_{i1}),\cdots,F_{d}(X_{id})\right),
\e*
which can be estimated by
\beg
\label{cross-validation}
\di \w {LSCV}({\mathbf N})= \int_{I^{d}} \left(\w{c}^{[{\mathbf N}]}({\mathbf u})\right)^2 d{\mathbf u} -\di \frac{2}{m}\di \sum_{i=1}^{m}\w{c}^{[{\mathbf N}]}_{-i}\left(\w F_{1}(X_{i1}),\cdots,\w F_{d}(X_{id})\right),
\en
yielding to the following estimator
of $\mathbf N$:
\be*
\di \w{\mathbf{N}}_{opt} =\argmin_{{\mathbf N}\in \N^d} \w{LSCV}({\mathbf N}), 
\e*
 where  $\di \w{c}^{[{\mathbf N}]}_{-i}$  is the leave-one-out copula density estimator without the data point $\mathbf{X}_i=(X_{i1},\cdots,X_{id})$.
Note that the expression (\ref{cross-validation})  has a similar form as the LSCV criterion used by 
\citet{bouezmarni2013bernstein}.

The proposition below gives an abridged form of $\w{LSCV}(\mathbf N)$  which is very useful to decrease its computation time in the numerical study.
\begin{proposition}\label{prop 5}
Fix $\m N\in\N^d$. We have  
\be*
\di \w{LSCV}(\mathbf N)= \frac{1}{n^2}\di\sum_{\mathbf m\leq \mathbf N}\left( \sum_{i=1}^{n}\prod_{j=1}^{d}Q_{m_j}^2\Big(\w F_{j}(X_{ij})\Big)-\di \frac{n+1}{n-1}\sum_{k\neq i}\prod_{j=1}^{d}Q_{m_j}\Big(\w F_{j}(X_{ij})\Big)Q_{m_j}\Big(\w F_{j}(X_{kj})\Big) \right).
\e*
\end{proposition}
The form of the $\w{LSCV}$ given in Proposition \ref{prop 5}  has the advantage to be easily evaluated numerically and it will be used in our simulation to select the value of $\mathbf N$ by minimising this expression.
\begin{proposition}\label{prop 6}
Fix $\m N\in\N^d$, independent of $n$. We have 
\be*
\di \E\left(\w{LSCV}(N)\right) &=&
    MISE(\w{c}^{[\mathbf{N}]})-\|c\|_{2}^{2}+\mathcal{O}(n^{-1}).
\e*
\end{proposition}
\begin{remark}
In the case where the margins are known
we have
\be*
\di \E\left(\w{LSCV}(N)\right) &=&
    MISE(\widetilde{c}^{[{\mathbf{N}}]})-\|c\|_{2}^{2}. 
\e*
\end{remark}

\subsection{Finite-sample performance of copula estimators}\label{monte carlo experiments}

To simplify our numerical study we 
fix $N_1=\cdots = N_d :=N$. It is therefore still possible to gain in precision by taking the time to choose a best combination of the components of the degree approximation $\m N$.

We use Monte Carlo simulations to demonstrate the potential of the copula estimator
$\di \w{C}^{[{\mathbf N}]}$ given in (\ref{estcop3})  where the optimal
parameter $\mathbf N$ is selected as described in Section \ref{subsection 3.3}.
This estimator will be denoted by CN.
We compare the performance with four competitors, namely:
\begin{itemize}
\item The empirical copula 
\citep{deheuvels1979fonction}, denoted Emp;
\item The empirical checkerboard copula 
\citep{carley2002new}, denoted Check;
\item The empirical Bernstein copula 
\citep{sancetta2004bernstein}, denoted Berns10  and Berns25 with smoothing parameter $k=10$ and $k=25$ respectively;
\item The empirical beta copula 
\citep{segers2017empirical}, denoted Beta.
\end{itemize}
We consider the classic copulas below in bivariate and trivariate dimension in the simulation series:
\begin{itemize}
\item  Clayton copula;
\item  Frank copula;
\item  Gaussian copula;
\item  Gumbel copula;
\item Independence copula;
\item Joe copula;
\item Student t-copula with degrees of freedom $\nu= 17$.
\end{itemize}
The readers may refer to  \citet{Nelsen} for the explicit functional forms and properties of these copulas.
We consider three levels of dependence with kendall's $\di \tau$'s equal to $\di\tau=0.3$ (low dependence),
$\di\tau=0.55$ (middle dependence) and $\di\tau=0.8$ (high dependence) for each copula model. We generate in both case
iid data using each copula model of sizes $\di n=500$ and $\di n=1000$.
In order to evaluate the quality of an estimator $\di \w{C}$ for a given copula $\di C$, we consider three performances measures:
the first one is the mean integrated absolute error (MIAE), defined by
\be*
\text{MIAE}(\w{C})=\E\left(\int_{I^d}\left\vert\w{C}^{[{\mathbf N}]}(\mathbf{u})-C(\mathbf{u})\right\vert d\mathbf{u} \right),
\e*
the second one is the mean integrated squared error (MISE), defined by
\be*
\text{MISE}(\w{C})=\E\left(\int_{I^d}\left\vert\w{C}^{[{\mathbf N}]}(\mathbf{u})-C(\mathbf{u})\right\vert^{2}d\mathbf{u}\right),
\e*
and the last one is the mean Kolmogorov-Smirnov error (MK-SE) defined by
\be*
\text{MK-SE}(\w{C})=\E\left(\sup_{\mathbf{u}\in I^d}\left\vert\w{C}^{[{\mathbf N}]}(\mathbf{u})-C(\mathbf{u})\right\vert\right),
\e*
estimated by the average over $M=1,000$ Monte Carlo replications of the approximate integration and Kolmogorov-Smirnov distance as follows
\be*
\di \Vert \w{C}^{[{\mathbf N}]}-C\Vert_{p}^{2}&\approx& \frac{1}{T^d}\sum_{j_1,\cdots,j_d=1}^{T-1}\left\vert \w{C}^{[{\mathbf N}]}({\mathbf j}/T)-C({\mathbf j}/T) \right\vert^{p},\, p=1,2,
\e*
\be*
\sup_{\mathbf{u}\in I^d}\left\vert\w{C}^{[{\mathbf N}]}(\mathbf{u})-C(\mathbf{u})\right\vert&\approx& \sup_{j_1,\cdots,j_d=1,\cdots,T-1}\left\vert\w{C}^{[{\mathbf N}]}({\mathbf j}/T)-C({\mathbf j}/T)\right\vert,
\e*
where ${\mathbf j}/T=(j_1/T,\cdots, j_d/T)$.

\paragraph{Two-dimensional case.}
In the  case where $d=2$, Tables \ref{MIAE-copula}
-\ref{MK-SE-copula} display the relatives
MIAE, MISE and MK-SE for the considered two sample sizes and
three levels of dependence.
It can be  observed that the copula
estimator $\w{C}^{[\mathbf{N}]}$ performs highly than the rest of
the estimators overall on these error criteria and for different
sample sizes and levels of dependence.
It is important to precise some comments on the individual comparisons:
\begin{itemize}
\item
For copulas with high dependence (Kendall's $\tau=0.8$), the
empirical beta copula can give better results  than $\w{C}^{[\mathbf N]}$, but with large standard
deviation in term of MIAE regardless of the  sample size.
In that case,  Bernstein copulas  with smoothing parameter $k=10$ and $k=25$  are really worse. We should try different  parameters $k$. Finally an automatic selection method is required to be able to detect the best Bernstein copula.

\item
For the independent copulas, the degree parameter
$N = 0$ is chosen all the time. In that case  $\w{C}^{[\mathbf N]}=\w{C}^{[\mathbf 0]}$ coincides with the true copula and it largely dominates
its competitors in all scenarios.
\item
Moreover, we can see  that our approach gives much better results
compared to other estimators with minimal standard deviation in
terms of relative MISE and MK-SE.
\item
We also remark that our method requires only small order of shifted
Legendre polynomial $\mathbf{N}$. This smoothing parameter
increase  with both the sample size ($n$)  and the level of
dependence ($\tau$).
\end{itemize}
\begin{center}
Table   \ref{MIAE-copula} here 
\\
Table \ref{MISE-copula} here 
\\
Table \ref{MK-SE-copula} here
\end{center}

\paragraph{Three dimensional case.} In the case where $d=3$, Tables
\ref{trivariate-copula200}  and 
\ref{trivariate-copula500} report the MIAE,MISE and MK-SE for   $n=200$ and $n=500$ and  for Kendall's
$tau=0.30$ and $0.8$, respectively. Our results are based on $M=100$ Monte Carlo
replications  and with  grid points  $j_i/T \in \{0.01,0.0712,\cdots,0.99\}$ for all $i=1,2,\cdots,d$.
Except for the Clayton copula where the beta estimator gave  the best result, the density estimator $\w{c}^{[\mathbf N]}$ dominates notably all the performances  with an extremely good results for the independence case since the choice $\mathbf N=0$ is almost always chosen and then the estimator fits exactly the density.

\begin{center}
 Table \ref{trivariate-copula200}  here
\\ Table \ref{trivariate-copula500} here 
\end{center}

\subsection{Finite-sample performance of copula density estimators}
Analogously to subsection \ref{monte carlo experiments}, we run Monte Carlo simulations to evaluate the performance
of the projection estimator of copula density $\di \w{c}^{[{\mathbf N}]}$ given by (\ref{copknown}). A finite-sample comparison with the following recent developments
in copula density estimators are considered.
\begin{itemize}
\item The probit-transformation estimator studied in 
\citet{Geenens}, denoted by $\w{c}_{Pt}^{1}$
for the local log-linear estimator and by $\w{c}_{Pt}^{2}$ for the local log-quadratic estimator;
\item The penalized hierarchical B-splines estimator studied in 
\citet{Kauermann}, denoted by  $\w{c}_{Ph}$;
\item The Bernstein copula density estimator studied in 
\citep{bouezmarni2013bernstein,bouezmarni2010asymptotic,janssen2014note}, denoted by $\w{c}_{Be10}$ and $\w{c}_{Be25}$
with smoothing parameter $k=10$ and $k=25$, respectively;
\item The thresholding estimators studied in 
\citet{Autin}, denoted by $\w{c}_{Th}^{~~l}$ for the local
thresholding estimator and by $\w{c}_{Th}^{~~b}$ for the block thresholding estimator. Haar wavelets are considered here;
\item The wavelet estimator studied in 
\citet{Genest}, denoted $\w{c}_{wa}^{}$. Haar wavelets are considered;
\item The beta kernel estimator studied in 
\citet{charpentier2007estimation},  denoted by $\w{c}_{Bk}$;
\item The Mirror reflection kernel estimator studied in 
\citet{gijbels1990estimating}, denoted by $\w{c}_{Mr}$.
\end{itemize}
The functions $\w{c}_{Pt}^{1}$, $\w{c}_{Pt}^{2}$, $\w{c}_{Ph}$,$\w{c}_{Be}$ and $\w{c}_{Bk}$ are provided in the R package \textbf{kdecopula} \citep{kdecop}.

Note that we are not looking at the thresholding and penalised hierarchical B-splines approaches due to the slowness
 of their calculation which makes them less competitive. 
We consider the same classic copulas and performance measures described above but narrowing to $M=100$
Monte Carlo replications and taking $j_i/T \in \{0.01,0.0712),\cdots,0.99\}$ for all $i=1,2,\cdots,d$.

Table \ref{MISE-density} reports the relative MISE and Table \ref{MK-SE-density} reports the relatives MK-SE for the considered two sample sizes and three levels of dependence.
They show that our approach clearly outperforms all the competitors considered for all sample sizes and levels of dependence. It appears that the degree parameter (truncation order) of our estimator has an influence on his performance. It increases when the level of dependence approaches one.  Note that the maximum optimal degree parameter, in all scenarios of the simulation study is $N=20$. But in the majority of cases this is between $N=0$ (for the independent case) and $N=10$.

\begin{center}
 Table  \ref{MISE-density} here 
 \\
 Table \ref{MK-SE-density} here 
\end{center}

\section{Real data applications}

\subsection{Insurance data}
\begin{figure}[!htbp]
\centering
\hspace*{-3cm}
\includegraphics[scale=0.7]{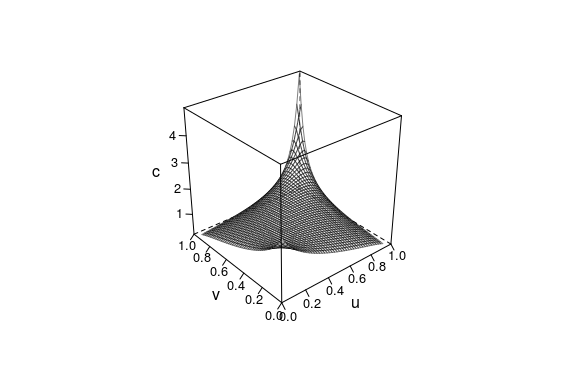}%
\hspace*{-3.5cm}
\includegraphics[scale=0.71]{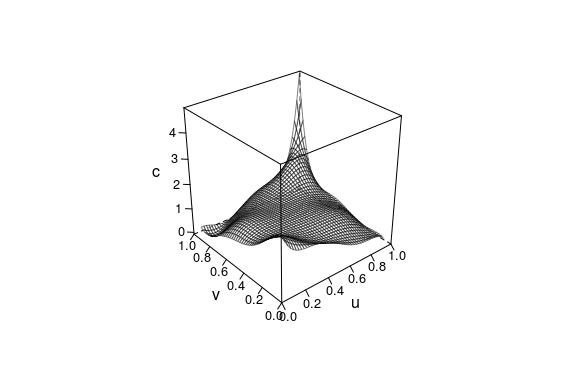}
\\[-1.25cm]
\hspace*{-2cm}
\includegraphics[scale=0.55]{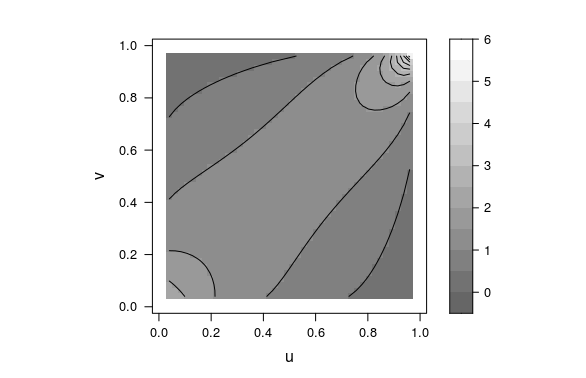}%
\hspace*{-1cm}
\includegraphics[scale=0.55]{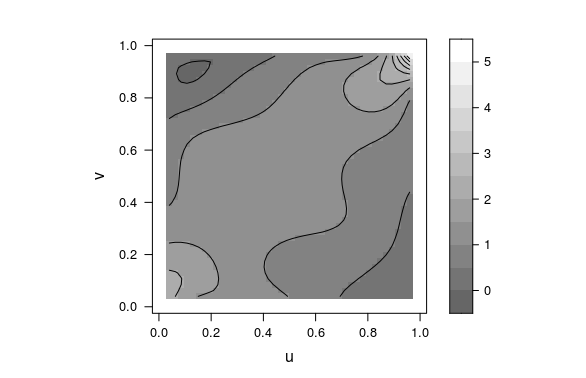}
\hspace*{-2cm}
\includegraphics[scale=0.5]{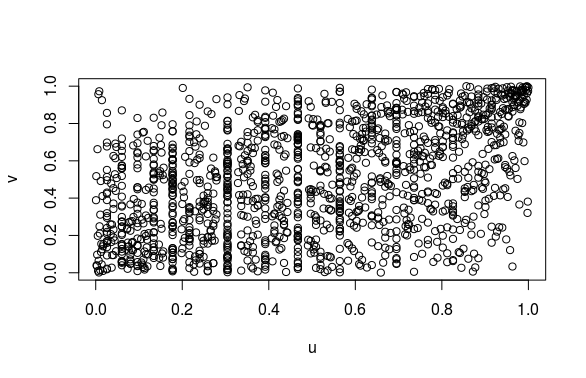}%
\hspace*{1cm}
\includegraphics[scale=0.4]{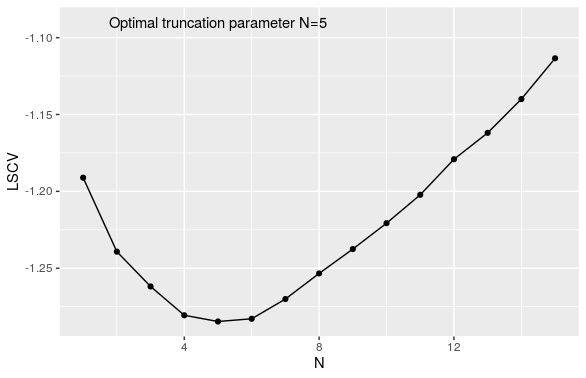}%
\\
\caption[]{Loss-ALAE data.  \\ Upper: Gumbel copula density with parameter $1.45$ (left);  $\w{c}^{[\mathbf{N}]}$ with $\mathbf{N}=5$ (right). \\
Middle: Gumbell  contours lines (left); $\w{c}^{[\mathbf{5}]}$ contours lines (right).\\  Bottom:  panel the rank-rank plot (left);
the LSCV function (right). }
\label{Fig:loss-alae}
\end{figure}

\begin{figure}[!htbp]
\includegraphics[scale=0.25]{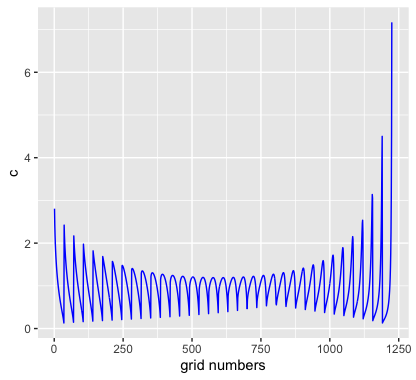}%
\includegraphics[scale=0.25]{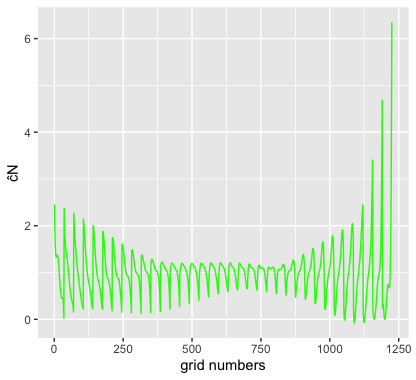}
\includegraphics[scale=0.25]{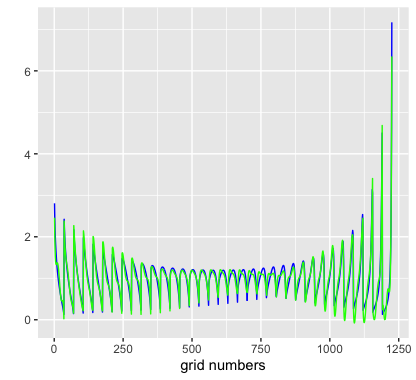}%
\caption[]{Loss-ALAE. 2-D line plot: Gumbel copula density with parameter $1.45$ (left); $\w{c}^{[5]}$ (middle); and their superposition (right).}
\label{dens alae point by point}
\end{figure}

\begin{figure}[!htbp]
\centering
\includegraphics[scale=0.25]{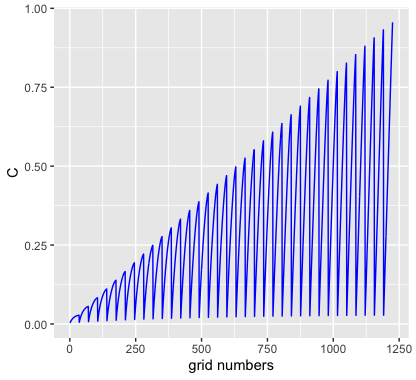}%
\includegraphics[scale=0.25]{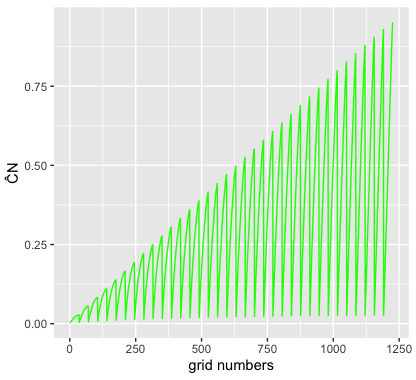}
\includegraphics[scale=0.25]{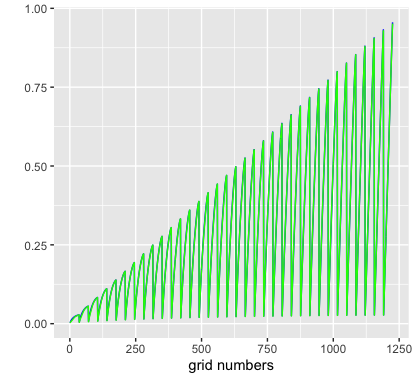}%
\caption[]{Loss-ALAE. 2-D line plot: Gumbel copula with parameter $1.45$ (left);  $\w{C}^{[5]}$ (middle);  and  their superposition (right).}
\label{copula alae point by point}
\end{figure}

We consider  a very  classical data set in the copula literature, namely the Loss-ALAE data set, that was collected by the US Insurance Services Office. It contains $1500$ general liability claims
each composed of the indemnity payment (Loss) and the allocated loss adjustment expense (ALAE).  We  excluded $34$ censored observations of the data set.  Many authors have used copulas
to model the dependence between the two variables in this data set, including 
\citet{frees1998understanding}, 
\citet{klugman1999fitting}, 
\citet{chen2005pseudo}, 
\citet{genest2006goodness}, 
\citet{denuit2006bivariate} and 
\citet{chen2010estimation}. The general conclusion is  that the Gumbel copula provides an adequate fit for these data.
Our purpose here is not to take back the analyses, but to confront our own estimator
to the  adjusted Gumbel copula. 
Figure \ref{Fig:loss-alae}  (right in the last row) shows the graph of the function LSCV using the selection rule prescribed in section \ref{subsection 3.3}.
The selected truncation parameter is $N=5$.
This figure also shows the proximity between our estimator $\w C^{[5]}$ and the Gumbell copula.  
To confirm the proximity  of  the copula as well as its density, 
we propose to compare them point by
point on a grid $\mathcal{G }$ with
$36 \times 36$ points uniformly chosen on $I^2$.  
The corresponding 1296 values associated to  the copula density estimator and to the theoretical Gumbel copula density are  represented in Figure \ref{dens alae point by point}.
The similarities between the two densities are 
clearly demonstrated, even at the  boundaries of $I^2$. 
The same analysis of  the copula is displayed in Figure
\ref{copula alae point by point}. It appears that the Gumbel copula and the estimator copula  merge
 very well.
Note that the grid numbers on the x-axis of these figures
represent the row number sequence in the data frame
$\mathcal{G}$. So each number corresponds to a couple
$(u_j, v_j)$.

\begin{center}
 Figure   \ref{Fig:loss-alae}  here 
 \\
 Figure \ref{dens alae point by point} here 
 \\
 Figure \ref{copula alae point by point} here 
\end{center}

In conclusion, even if (\ref{condition}) is not  satisfied  for a Gumbel copula density we obtain a very efficient estimators. As explained in Appendix \ref{appendL2} this is due to the fact that for any $\epsilon >0$,  writing $c$ the Gumbel  copula density, there exists a function $ c' \in \mathbb{L}^2([0,1]^d)$ such that, for all $\m u$, $|c(u)- c'(u)|<\epsilon$. Their  associated copula coefficients  can also be as close as desired. Then our procedure yields a good estimators of $ c'$ which is also a good approximation of $c$.

\subsection{Financial data}

\begin{figure}[!htbp]
\centering
\hspace*{-3cm}
\includegraphics[scale=0.35]{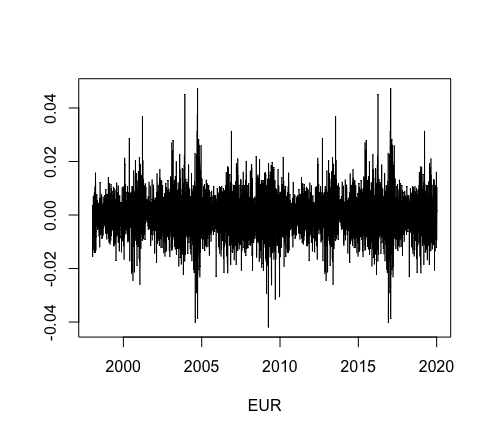}%
\hspace*{-0.5cm}
\includegraphics[scale=0.35]{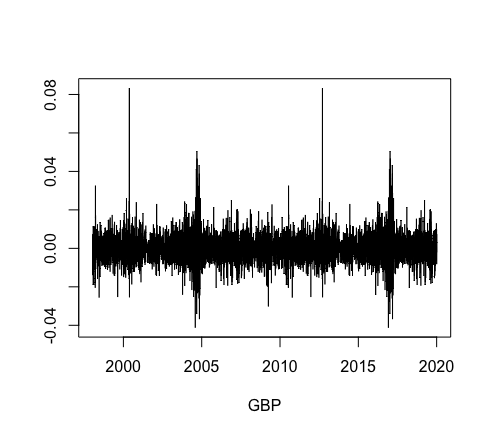}
\\[-0.5cm]
\hspace*{-3cm}
\includegraphics[scale=0.35]{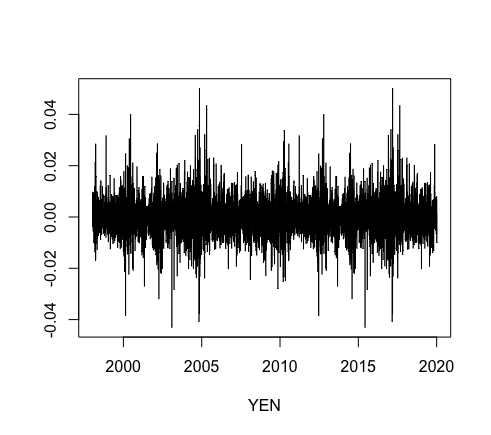}%
\hspace*{-0.5cm}
\includegraphics[scale=0.35]{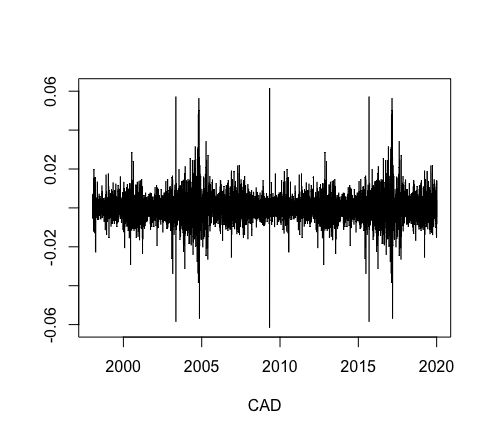}
\\[-0.5cm]
\hspace*{-3cm}
\includegraphics[scale=0.35]{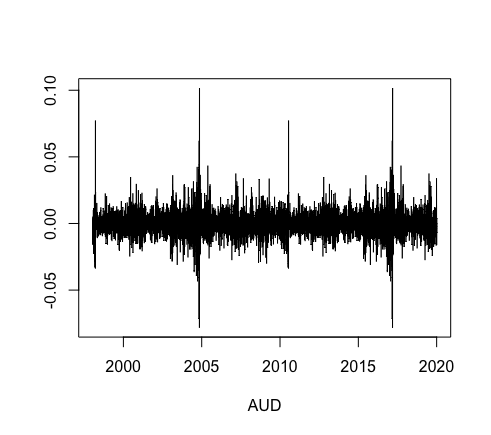}%
\hspace*{-0.5cm}
\includegraphics[scale=0.35]{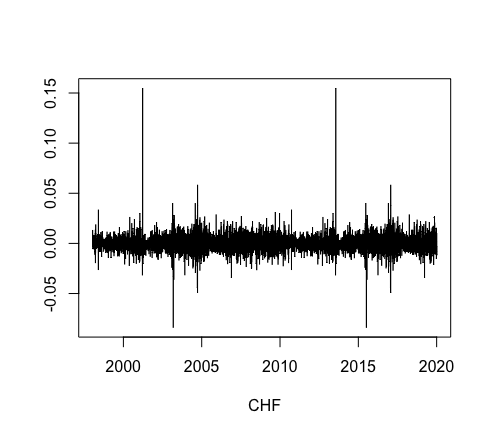}%
\\[-0.3cm]
\caption[]{Time plots of daily log returns of Euro, Great British pound, Japanese yen, Canadian Dollar, Australian Dollar and Swiss franc}
\label{Fig:currency}
\end{figure}

In this data analysis, we study the dependence structure
between the series of the most used exchange rates in foreign
exchange markets. Data are available from IMF (International Monetary Fund) rates database and
consists in daily currency exchange rates from January
$1^{st}$, $1994$ to July $31^{th}$, $2020$ for a total of
$5551$ business days (any day except weekends and some
holidays). We exclude from our analysis $1049$ observations
with at least one missing value. By the standard continuously
compounded return formula, we consider the log-returns of six
exchanges rates: the Euro (EUR), the Great British pound
(GBP), the Japanese yen (YEN), the Canadian dollar (CAD), the
Swiss franc (CHF) and the Australian dollar (AUD).
The time plots of these returns are shown in Figure
\ref{Fig:currency}. 
These time plots of log returns have a similarly behave with
time and show the stylized fact of clustering volatility.
We  also notice the appearance of extreme values.
Summary statistics for these returns are displayed in Table
\ref{sum_data}. It reveals that the mean and median of
log-returns of these six exchange rates are very close to
zero. Their distribution are positively skewed (right tailed)
and leptokurtic (kurtosis value higher than the kurtosis of
normal distribution whose value is three). The CHF has the
highlest kurtosis and skewness, the EUR the lowest kurtosis and the CAD the lowest skewness.
As it was expected due to more frequent occurrence of extreme
values in Figure \ref{Fig:currency}, 
none of daily log returns
passed the Jarque-Bera (JB) test ($P_{value}<0.001$) where the
null hypothesis is that the log-return series are normally
distributed.
As a result, the use of standard models is not appropriate in this context because one of the most important assumptions of several models of the financial time series is the normality of returns.

Tables \ref{parameter_data} and   \ref{kendall_data} give
respectively the optimal degree parameters and Kendall's
rank correlations between the six daily log-returns.
It shows that the proposed estimators don't need a high
optimal truncation order. The largest degree 
parameter is $N=16$ (EUR/CHF) and the smallest one is $N=3$ (EUR/YEN).
Futher, various Kendall correlations are small values, positive for some pairs of daily log returns and negative for others.
The pair EUR/CHF has the highest Kendall correlation in absolute value ($\tau=0.6308960$). We note that a  large value for $N$ increases the bias and consequently decreases the asymmetrical structure.

\begin{figure}[!htbp]
\centering
\vspace*{-3cm}
\hspace*{-3cm}
\includegraphics[scale=0.55]{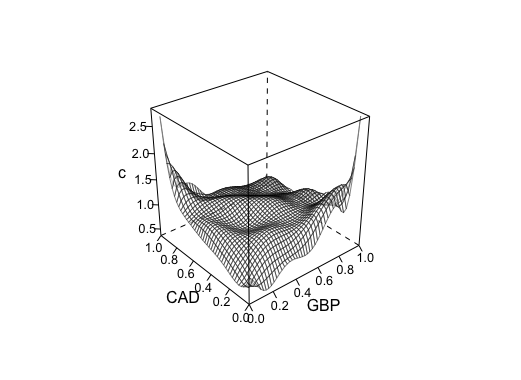}%
\hspace*{-0.5cm}
\includegraphics[scale=0.55]{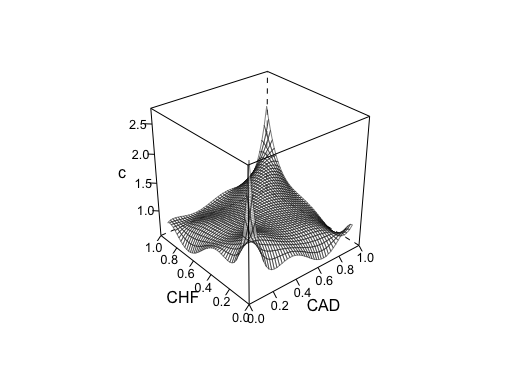}
\\[-0.5cm]
\hspace*{-3cm}
\includegraphics[scale=0.55]{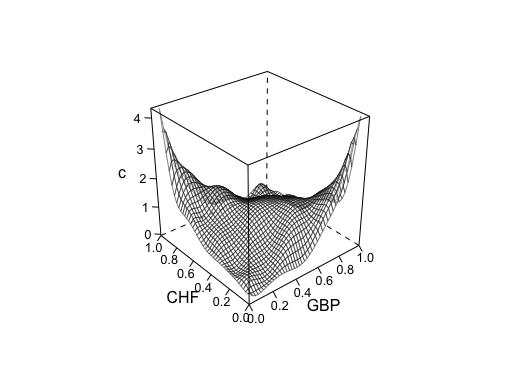}%
\hspace*{-0.5cm}
\includegraphics[scale=0.55]{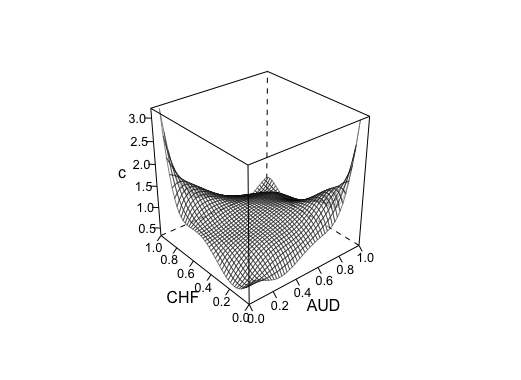}
\\[-0.5cm]
\hspace*{-3cm}
\includegraphics[scale=0.55]{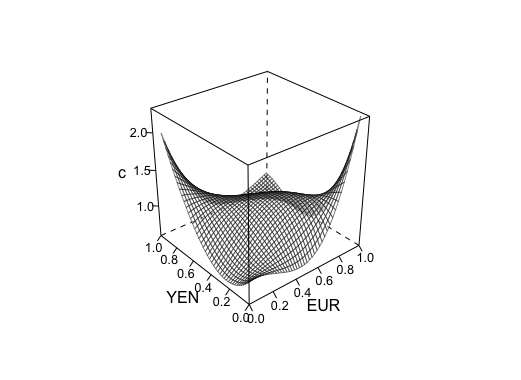}%
\hspace*{-0.5cm}
\includegraphics[scale=0.55]{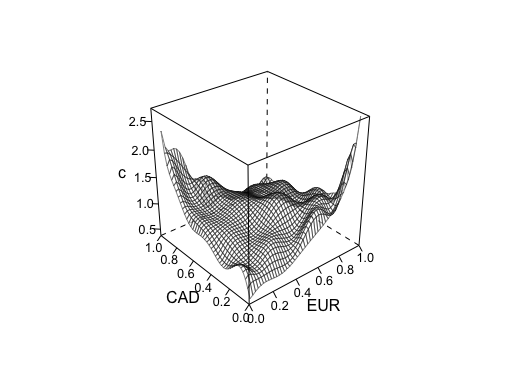}%
\\[-0.3cm]
\caption[]{Estimators of the copula density for Canadian Dollar/Great British pound, Swiss franc/Canadian Dollar, Swiss franc/Great British pound, Swiss franc/Australian Dollar, Japanese yen/Euro and  Canadian Dollar/Euro }
\label{Fig:dens-currency1}
\end{figure}

\begin{figure}[!htbp]
\centering
\vspace*{-3cm}
\hspace*{-3cm}
\includegraphics[scale=0.55]{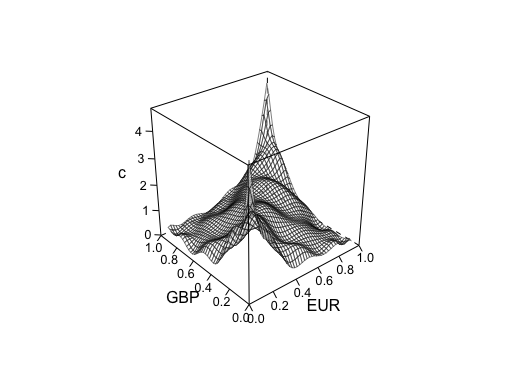}%
\hspace*{-0.5cm}
\includegraphics[scale=0.55]{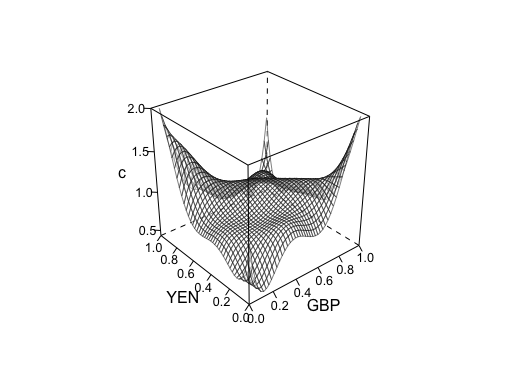}
\\[-0.5cm]
\hspace*{-3cm}
\includegraphics[scale=0.55]{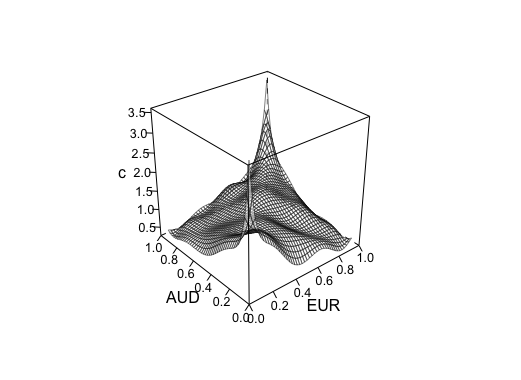}%
\hspace*{-0.5cm}
\includegraphics[scale=0.55]{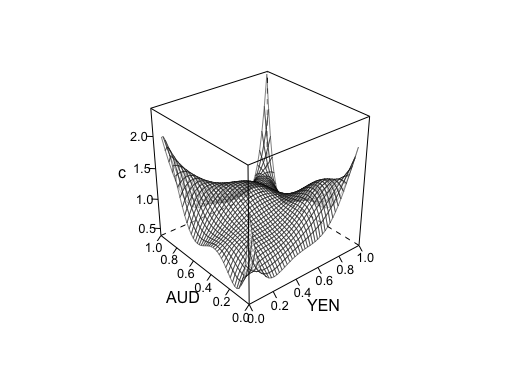}
\\[-0.5cm]
\hspace*{-3cm}
\includegraphics[scale=0.55]{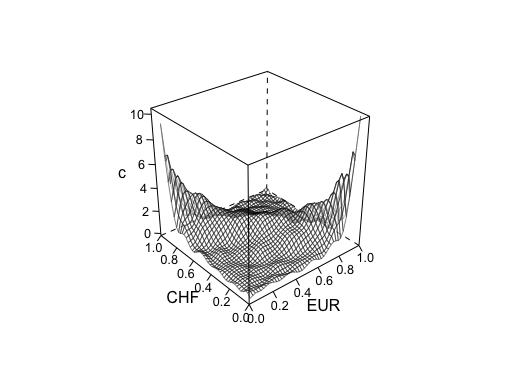}%
\hspace*{-0.5cm}
\includegraphics[scale=0.55]{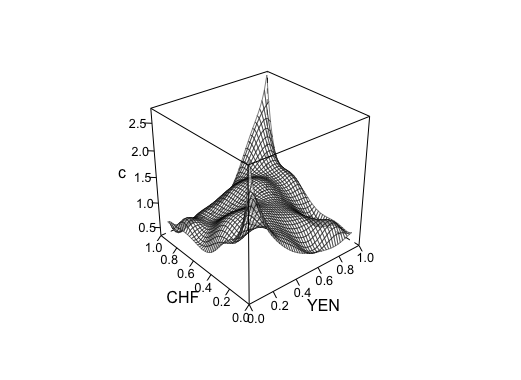}%
\\[-0.3cm]
\caption[]{Estimators of the copula density for Euro/Great British pound, Japanese yen/Great British pound, Australian Dollar/Euro, Australian Dollar/Japanese yen,  Swiss franc/Euro and  Swiss franc/Japanese yen}
\label{Fig:dens-currency2}
\end{figure}

As it was expected from the simulation studies in Section \ref{simul-studies}, we notice a close link between the Kendall coefficients and the optimal degree parameters of proposed estimators. Both indicators increase or decrease simultaneously.
These information are useful for the development of risk diversification strategies since the investment in a portfolio of various assets reduces risk, particularly in exchange rate management. The study can be extended to calculations of the
tail value at risk and risk measure expected shortfall.
The different pairs daily log return are presented in Figure \ref{Fig:dens-currency1} and \ref{Fig:dens-currency2}. We see that the proposed approach captures very well the asymmetry in the dependence structure between these pairs of daily log returns exchange rate.

\begin{center}
Table   \ref{sum_data} here 
\\
Table \ref{parameter_data}  here 
\\
Table  \ref{kendall_data} here
\end{center}

\begin{center}
Figure \ref{Fig:currency} here 
\\ 
Figure  \ref{Fig:dens-currency1}  here 
 \\
 Figure \ref{Fig:dens-currency2} here 
\end{center}

\section{Conclusion}
In this paper a very simple nonparametric estimator of the copula density based on shifted Legendre polynomials is proposed.  A new copula estimator is then deduced based on copula coefficients.  
Both estimators are easy to implement with an automatic selection of degree approximation. A R program is available on \href{https://github.com/yvesngounou/}{Github-yvesngounou}.
Various theoretical properties are demonstrated, including asymptotic convergence and optimality in the minimax 
sense. 
Experimental results clearly demonstrated the performance of the proposed estimators. 
The proposed copula estimators seem  to outperform 
all its recent competitors found in the literature  
according to the MIAE, MISE and MK-SE criterion in various scenarios. These experiment results also showed the superiority of the proposed copula density estimator. 
These performances could also  be improved when considering different orders for  $N_1,N_2,\cdots,N_d$ instead of equal values $N_1=\cdots = N_d$. This seems  possible at a reasonable  computationally  cost. 
In the independence case, the estimator is extremely accurate because the correct copula is generally selected. This result may lead to thinking about the construction of an independence test. Moreover, the correspondence  that we have shown between the  copula coefficients and the Spearman's rho is an interesting direction for future research.
Real situations in actuarial and financial frameworks have demonstrated the adaptability of the proposed method.
Finally, this approach by orthogonal projections  should also allow us to construct test statistics   to compare copulas.

We close this conclusion by  recalling  that this approach is nonparametric and then hypothesis (\ref{condition})  cannot be verified. However, the method works extremely well even for cases where the density is not square integrable which makes it very useful in practice. We explain this phenomenon  in Appendix \ref{appendL2}. 




\bibliographystyle{apalike}

\bibliography{papercopula}




%




\section*{APPENDIX: Proofs}



We will denote 
\be*
\di \sup \lvert \w F_j -F_j\rvert:= \lVert \w F_j -F_j\rVert_{\infty}, \quad \pi(\bm m)=:\di \prod_{j=1}^{d} m_j, \quad  \max\{N_j; j:=1,\cdots,d\}:=N_{max}.
\e*
Let begin with the following useful lemmas

\begin{lemma}\label{lem_Op}
\be*
 \di \lVert \w F_j -F_j\rVert_{\infty}=o_{\P}(1) \text{ and } 
 \di \lVert \w F_j -F_j\rVert_{\infty}=\mathcal{O}_{\P}(\sqrt{n^{-1}})
\e*
\end{lemma}

According to the Massart inequality (\cite{massart1990tight}), we get
\beg \label{massart}
 \di \forall\epsilon>0, \quad \P\left( \lVert \w F_j -F_j\rVert_{\infty}>\epsilon\right)\leq 2e^{-2n\epsilon^2}
\en
 and the results follow.
 \begin{lemma}\label{lem_exp_Op}
 \be*
\di  \E \lVert \w F_j -F_j\rVert_{\infty}^2 \leq n^{-1} \quad \text{ and }\quad
\di  \E\lVert \w F_j -F_j\rVert_{\infty}\leq \sqrt{n^{-1}} 
 \e* 
 \end{lemma}
 Denoting $Y:=\lVert \w F_j -F_j\rVert_{\infty}$ and using (\ref{massart}), we see that
 \be*
\di  \E (Y)\leq \sqrt{ \E (Y^2)}=\sqrt{\int_{0}^{+\infty}\P( Y>\sqrt{s}) ds}\leq \sqrt{\int_{0}^{+\infty}2e^{-2ns^2} ds}=\sqrt{n^{-1}} 
\e* 
 and the results follow.
 \begin{lemma}\label{lem_majora}
 For all $\alpha\geq 1$,
 \be* 
 \di \sum_{{\mathbf m}\leq {\mathbf N}}\big(\pi(\bm m)\big)^{\alpha}\leq N_{max}^{d(\alpha+1)}
 \e*
 \end{lemma}
We have 
 \be* 
 \di \sum_{{\mathbf m}\leq {\mathbf N}}\big(\pi(\bm m)\big)^{\alpha}= 
 \prod_{j=1}^{d} \sum_{m_j=0}^{N_j} m^{\alpha-1}_j m_j\leq N_{max}^{d(\alpha-1)}\prod_{j=1}^{d} \sum_{m_j=0}^{N_j} m_j\leq N_{max}^{d(\alpha-1)}\prod_{j=1}^{d} \frac{N_j(N_j+1)}{2} \leq N_{max}^{d(\alpha-1)}\prod_{j=1}^{d}N_j^2\leq N_{max}^{d(\alpha+1)}.
 \e*

 \begin{lemma}\label{lem_legendre}
For all $u \in I$ we have 
\beg
\label{inegLegendre}
|Q_m(u)| \leq \eta_1 m^{1/2}
&{\rm and } & |Q'_m(u) | \leq \eta_2 m^{5/2}\qquad \forall m>0, 
\en
where $\eta_1= \sqrt{3}$ and $\eta_2= 2\sqrt{3}$.


\end{lemma}
 See for instance inequalities $22.14$ on page $791$ of \citet{abramowitz1970handbook} or \citet{boas1969inequalities}.
 .  

\subsection*{Proof of Proposition \ref{prop0}}

From (\ref{densite2}) we have $h_{\mathbf U}=c$. Since $(\mathbf{Q}_{\mathbf m})_{\mathbf m \in \N^d}$ forms a dense orthogonal basis with respect to the uniform measure on $I^d$, by (\ref{condition}) we have
\begin{align*}
h_{\mathbf U}(u_1,\cdots, u_d) & =
\di\sum_{\mathbf m \in \N^d}
\di\left(\int c(v_1,\cdots, v_d)\mathbf{Q}_{\mathbf m}(\mathbf v)  dv_1 \cdots dv_d \right)\mathbf{Q}_{\mathbf m}(\mathbf u),
\end{align*}
which gives (\ref{classic}), and
(\ref{cop2}) follows by integration.
\subsection{Proof of Proposition \ref{bona_fide}}
Assume that $u_i=0$. Then we have
\be*
\di \w{C}^{[{\mathbf N}]}({\mathbf u})=
\sum_{{\mathbf m}\leq {\mathbf N}}\w{\rho}_{{\mathbf m}}\left(\int_{0}^{0}Q_{m_i}(x_i)\mu(dx_{i})\right)\prod_{\substack{j=1 \\j\neq i}}^{d}\int_{0}^{u_j}Q_{m_j}(x_j)\mu(dx_{j})=0,
\e*
and $i)$ is proved.
%
To prove $ii)$ assume that $\di \mathbf{u}=(1, \cdots,1,u_i,1, \cdots,1)$. Then
\be*
\di \w{C}^{[{\mathbf N}]}({\mathbf u})&=&
\sum_{{\mathbf m}\leq {\mathbf N}}\w{\rho}_{{\mathbf m}}\left(\int_{0}^{u_i}Q_{m_i}(x_i)\mu(dx_{i})\right)\prod_{\substack{j=1 \\j\neq i}}^{d}\int_{0}^{1}Q_{m_j}(x_j)\mu(dx_{j})\\
&=&\sum_{{\mathbf m}\leq {\mathbf N}}\w{\rho}_{{\mathbf m}}\left(\int_{0}^{u_i}Q_{m_i}(x_i)\mu(dx_{i})\right)\prod_{\substack{j=1 \\j\neq i}}^{d}\delta_{0,m_{j}}\\
&=&\w{\rho}_{\mathbf 0} \int_{0}^{u_i}Q_{0}(x_i)\mu(dx_{i})=
\int_{0}^{u_i}1\;\mu(dx_{i})=u_i.
\e*
$iii)$ is immediate from  $\w{\rho}_{\mathbf 0}=1$ and $Q_{0}=1$.

\subsection*{Proof of Proposition \ref{pro_uniform}}
By orthogonality of the polynomials we have
\be*
\di\int_{I^{d}}\w{c}^{[{\mathbf N}]}({\mathbf u})\mathbf {\mu}(d{\mathbf u})&=&
\sum_{{\mathbf m}\leq {\mathbf N}}\w{\rho}_{{\mathbf m}}\prod_{j=1}^{d}\int_{0}^{1}Q_{m_j}(u_j)\mu(du_{j})\\
&=&\sum_{{\mathbf m}\leq {\mathbf N}}\w{\rho}_{{\mathbf m}}\prod_{j=1}^{d}\delta_{0,m_{j}}\\
&=& 1
\e*
which yields $i)$.
The proof of $ii)$ is very similar since $\rho_{\mathbf m}=0$ if exactly  $d-1$ components of $\mathbf m$ are null.  $iii)$ is immediate.

\subsection*{Proof of Proposition \ref{theo spearman}}

We have
\be*
\di \rho_{11}&=&\E\left(Q_{1}(F_{X_1}(X_1))Q_{1}(F_{X_2}(X_2)) \right)\\
&=& \E_{C}\left(Q_{1}(U_1))Q_{1}(U_2) \right)\\
&=& \int_{0}^{1}\int_{0}^{1}\left(\sqrt{3}(2u_1-1).\sqrt{3}(2u_2-1)\right)C(du_1,du_2)\\
&=&3 \int_{0}^{1}\int_{0}^{1}\left(4u_{1}u_{2}-2u_{1}-2u_{2}+1\right)C(du_1,du_2)\\
&=&3 \left(\int_{0}^{1}\int_{0}^{1}4u_{1}u_{2}dC(u_1,u_2)-2(\frac{1}{2})-2(\frac{1}{2})+1\right)\\
&=&12\int_{0}^{1}\int_{0}^{1}4u_{1}u_{2}dC(u_1,u_2)-3\\
&=& \rho_C.
\e*

\subsection*{Proof of Proposition \ref{prop 1}}
We have
\be*
\di \w{c}^{[{\mathbf N}]}({\mathbf u})-\tilde{c}^{[{\mathbf N}]}({\mathbf u}) &=&\sum_{{\mathbf m}\leq {\mathbf N}}(\w{\rho}_{{\mathbf m}}-
\tilde{\rho}_{{\mathbf m}}\mathbf{Q}_{\mathbf{m}}({\mathbf u}))\\
&=&\di \sum_{{\mathbf m}\leq {\mathbf N}}\frac{1}{n}\sum_{i=1}^{n}\Bigg(\prod_{j=1}^{d}Q_{m_j}(\w{F}_{j}(X_{ij}))-\prod_{j=1}^{d}Q_{m_j}(F_{j}(X_{ij}))\Bigg)\mathbf{Q}_{\mathbf{m}}({\mathbf u}).
\e*

 By Taylor expansion, we obtain
\beg
|\prod_{j=1}^{d}Q_{m_j}(\w{F}_{j}(X_{ij}))-\prod_{j=1}^{d}Q_{m_j}(F_{j}(X_{ij}))|
&\leq&  \sum_{j=1}^{d} \eta_2 m_j^{5/2} \di\prod_{i\neq j} \eta_1 m_i^{1/2} \|\w{F}_{j}-F_{j}\|_{\infty}
\nonumber
\\
& \leq &
M({\mathbf m}) \sum_{j=1}^{d} \|\w{F}_{j}-F_{j}\|_{\infty},
\label{Taylor}
\en
where $M({\mathbf m}) = \eta_2\eta_1^{d-1}\pi({\mathbf m})^{1/2}\max(m_j)^{2}$. 
It follows that
\beg \label{eq1prop1}\nonumber
\di |\w{c}^{[{\mathbf N}]}({\mathbf u})-\tilde{c}^{[{\mathbf N}]}({\mathbf u})|
&\leq& \di \sum_{{\mathbf m}\leq {\mathbf N}}M({\mathbf m}) \max|\mathbf{Q}_{{\mathbf m}}|\sum_{j=1}^{d} \|\w{F}_{j}-F_{j}\|_{\infty}
\\\nonumber
&\leq& \di  \eta_2\eta_1^{2d-1} N_{max}^2 \sum_{{\mathbf m}\leq {\mathbf N}}\pi({\mathbf m})\sum_{j=1}^{d} \|\w{F}_{j}-F_{j}\|_{\infty}
\\
&\leq&  \eta_2\eta_1^{2d-1} N_{max}^{2d+2} \sum_{j=1}^{d} \|\w{F}_{j}-F_{j}\|_{\infty}, 
\en
from Lemma \ref{lem_majora}, which yields $i)$  since $\|\w{F}_{j}-F_{j}\|_{\infty} =  \mathcal{O}_{\P}\big(\sqrt{n^{-1}}\big)$ according to Lemma \ref{lem_Op}.
\\
We obtain $ii)$ analogously.

\subsection*{Proof of Proposition \ref{prop 4}.}

By orthogonality of the Legendre polynomials we have
\be*
\| \w{c}^{[\mathbf{N}]}-c\|_{2}^{2}
&=&\sum_{{\mathbf m}\leq {\mathbf N}}\left(\w{\rho}_{{\mathbf m}}-\rho_{{\mathbf m}}\right)^2
+
\sum_{{\mathbf m}\not\leq {\mathbf N}}\rho_{{\mathbf m}}^2,
\e*
which gives result.

\subsection*{Proof of Corollary \ref{relation_mise}}
By orthogonality, we have
\beg
\label{eqMISE0}
MISE(\w{c}^{[\mathbf{N}]})- MISE(\tilde{c}^{[\mathbf{N}]})&=&\E\left(\sum_{\mathbf{m} \leq \mathbf{N}}(\w{\rho}_{\mathbf{m}}-\rho_{\mathbf{m}})^2\right)
-\E\left(\sum_{\mathbf{m} \leq \mathbf{N}}(\tilde{\rho}_{\mathbf{m}}-\rho_{\mathbf{m}})^2\right)
\nonumber
\\ 
&=&\E\sum_{\mathbf{m} \leq \mathbf{N}}(\w{\rho}_{\mathbf{m}}-\tilde{\rho}_{\mathbf{m}})
\left((\w{\rho}_{\mathbf{m}}-\tilde{\rho}_{\mathbf{m}})+2(\tilde{\rho}_{\mathbf{m}}-\rho_{\mathbf{m}})\right).
\en
From (\ref{Taylor}), we first have
\be*
|\w{\rho}_{\mathbf{m}}-\tilde{\rho}_{\mathbf{m}}|
& \leq &
M({\mathbf m}) \sum_{j=1}^{d} \|\w{F}_{j}-F_{j}\|_{\infty},
\e*
and since  $$ \|\w{F}_{j}-F_{j}\|_{\infty} =  \mathcal{O}_{\P}\big(\sqrt{n^{-1}}\big)$$ it follows that
\beg \label{cor2-asum1}
\w{\rho}_{\mathbf{m}}-\tilde{\rho}_{\mathbf{m}}=\mathcal{O}_{\P}(\sqrt{n^{-1}}).
\en
Furthermore, by the Central Limit Theorem we have
\beg \label{cor2-asum2}\nonumber
&&\sqrt{n}\left(\tilde{\rho}_{\mathbf{m}}-\rho_{\mathbf{m}}\right)
\xrightarrow{\enskip \mathcal{L}\enskip} \mathcal{N}\left(0,\mathbb{V}\big(\prod_{j=1}^{d}Q_{m_j}
(F_{j}(X_{ij}) \big)\right),
\en
where $\di \mathbb{V}\big(\prod_{j=1}^{d}Q_{m_j}
(F_{j}(X_{ij}) \big)<\infty$ since $\bm{N}$ is fixed and independent of $n$. 
Then
\beg
\label{eqOP}
\tilde{\rho}_{\mathbf{m}}-\rho_{\mathbf{m}}
&=&\mathcal{O}_{\P}\big(\sqrt{n^{-1}}\big).
\en
Combining (\ref{eqMISE0}), (\ref{Taylor}) and (\ref{eqOP}) we get
\be*
MISE(\w{c}^{[\mathbf{N}]})= MISE(\tilde{c}^{[\mathbf{N}]})+
\mathcal{O}\big(n^{-1}\big).
\e*

\subsection*{Proof of Corollary \ref{corol 1}}
From Proposition \ref{prop 4}, we have
\be*
\E\left( \| \w{c}^{[\mathbf{N}]}-c\|_{2}^2\right) =\E\left(\sum_{\mathbf{m} \leq \mathbf{N}}(\w{\rho}_{\mathbf{m}}-\rho_{\mathbf{m}})^2\right) +  \di \sum_{\mathbf{m} \not\leq \mathbf{N}}\rho_{\mathbf{m}}^2.
\e*
From the elementary inequality $\di 2ab \leq a^2 + b^2$ for
any $a,b \in \R$, we get
\beg
\label{inegcor1}
\E\left( \| \w{c}^{[\mathbf{N}]}-c\|_{2}^2\right) &\leq &2\E\left(\sum_{\mathbf{m} \leq \mathbf{N}}(\w{\rho}_{\mathbf{m}}- \tilde{\rho}_{\mathbf m} )^2\right) +2\E\left(\sum_{\mathbf{m} \leq \mathbf{N}}(\tilde{\rho}_{\mathbf{m}}-
\rho_{\mathbf{m}})^2\right) +  \di \sum_{\mathbf{m} \not\leq \mathbf{N}}\rho_{\mathbf{m}}^2
\nonumber \\
& := &
2A_1 + 2A_2 + A_3.
\en
Clearly, $A_3 \rightarrow 0$, as $N \rightarrow \infty$ as a tail of a convergent series.


Moreover by (\ref{inegLegendre}) we have
\begin{align}
\nonumber
 A_2 & =
\sum_{\mathbf{m} \leq \mathbf{N}}
 \frac{1}{n} \V(\di\prod_{j=1}^dQ_{m_j}(F_j(X_j)))
 \\\nonumber
 & \leq
\sum_{\mathbf{m} \leq \mathbf{N}} \di \frac{1}{n}
 \eta_{1}^{2d} \di\prod_{j=1}^d m_j
 \\\label{eq1}
 & \leq
 \di\frac{\eta_{1}^{2d} N_{max}^{2d}}{n},
\end{align}
By {\bf (H)}  $A_2 \rightarrow 0$ as $n \rightarrow \infty$. 
To conclude, we combine (\ref{Taylor})  with (\ref{inegLegendre}) to obtain
\begin{align*}
A_1 & =
\sum_{\mathbf{m} \leq \mathbf{N}}\di\frac{1}{n^2}
\di\sum_{i=1}^n \di\sum_{i'=1}^n
\E\Bigg[\Bigg(\prod_{j=1}^{d}Q_{m_j}(\w{F}_{j}(X_{ij}))-\prod_{j=1}^{d}Q_{m_j}(F_{j}(X_{ij}))\Bigg)
\\
& \ \ \ \ \ \ \ \ \ \ \times \Bigg(\prod_{j=1}^{d}Q_{m_j}(\w{F}_{j}(X_{i'j}))-\prod_{j=1}^{d}Q_{m_j}(F_{j}(X_{i'j})))\Bigg)\Bigg]
\\
& \leq
\sum_{\mathbf{m} \leq \mathbf{N}}\di\frac{1}{n^2}
\di\sum_{i=1}^n \di\sum_{i'=1}^n
\E\Bigg(
M({\mathbf m})^2 \left(\sum_{j=1}^{d} \|\w{F}_{j}-F_{j}\|_{\infty}\right)^2\Bigg)
\\
& \leq
d\sum_{\mathbf{m} \leq \mathbf{N}}M({\mathbf m})^2\sum_{j=1}^{d}\E
 \left( \|\w{F}_{j}-F_{j}\|^{2}_{\infty}\right), 
 \end{align*}
 according to Hölder's sum inequality.

Since $M({\mathbf m}) = \eta_1^{d-1}\eta_2 \pi({\mathbf m})^{1/2}\max(m_j)^{2}$ and $\E
 \left( \|\w{F}_{j}-F_{j}\|^{2}_{\infty}\right) \leq n^{-1}$,  we see from to Lemma \ref{lem_exp_Op} that there exists a constant $\eta>0$ such that
\begin{align}
\nonumber
A_1 & \leq \frac{d^2 \eta_1^{2d-2}\eta^2_2}{n} \sum_{\mathbf{m} \leq \mathbf{N}}\pi(\bm m)\max(m_j)^{4}
\\\nonumber
& \leq \frac{d^2 \eta_1^{2d-2}\eta^2_2 N_{max}^4}{n} \sum_{\mathbf{m} \leq \mathbf{N}}\pi(\bm m)
\\ \label{eq2}
 & \leq \frac{d^2 \eta_1^{2d-2}\eta^2_2 N_{max}^{2d+4}}{n} ,
\end{align}
and then $A_1 \rightarrow 0$ as $n \rightarrow \infty$.

\subsection*{Proof of Proposition \ref{prop 2}}
We have 
\be*
\di\E\left(\w{c}^{[{\mathbf N}]}({\mathbf u})\right)=\E \left(\tilde{c}^{[{\mathbf N}]}({\mathbf u})\right) +
\di\E\left(\w{c}^{[{\mathbf N}]}({\mathbf u})-\E\tilde{c}^{[{\mathbf N}]}({\mathbf u})\right). 
\e*
Since $
\E\left(\tilde{\rho}_{{\mathbf m}}\right) =\rho_{{\mathbf m}}$, we have
\be*
\di \E\left(\tilde{c}^{[{\mathbf N}]}(\mathbf{u})\right) &=& \sum_{{\mathbf m}\leq {\mathbf N}}\E\left(\tilde{\rho}_{{\mathbf m}}\right)\prod_{j=1}^{d}Q_{m_j}(u_j)=\sum_{{\mathbf m}\leq {\mathbf N}}\rho_{{\mathbf m}}\prod_{j=1}^{d}Q_{m_j}(u_j),
\e*
which tends to $c\left(\mathbf{u}\right)$  as $n$ tends to  $\infty$.
From (\ref{eq1prop1}), \text{ according to Lemma }\ref{lem_exp_Op},  we have
\be*
\di |\E\w{c}^{[{\mathbf N}]}({\mathbf u})-\E\tilde{c}^{[{\mathbf N}]}({\mathbf u})|
&\leq& \di \eta_2\eta_1^{2d-1} N_{max}^{2d+2} \sum_{j=1}^{d}\E \|\w{F}_{j}-F_{j}\|_{\infty}
\\\nonumber
&\leq& d \di  \eta_2\eta_1^{2d-1}\frac{ N_{max}^{2d+2}}{\sqrt{n}} .
\e*
Under assumption $(\bm{H})$, we deduce that 

\be*
\di \E\left(\w{c}^{[{\mathbf N}]}(\mathbf{u})\right) \longrightarrow  c(\mathbf{u}) 
 \text{ as } n \to \infty.
\e*
and analogously
\be*
\di  \E\left(\w{C}^{[{\mathbf N}]}(\mathbf{u})\right) \longrightarrow C(\mathbf{u})
 \text{ as } n \to \infty.
\e*

\subsection*{Proof of Proposition  \ref{theor1}}


Combining (\ref{inegcor1}) with (\ref{eq1})- (\ref{eq2}) yields  
%
%

%

\be*
\nonumber 
\di \E (\| \w{c}^{[{\mathbf{N}}]}-c\|_2^2 )
&\leq &
\di 2\frac{d^2 \eta_1^{2d-2}\eta^2_2 N_{max}^{2d+4}}{n}  + 2\di\frac{\eta_{1}^{2d} N_{max}^{2d}}{n} + 
\di \sum_{\mathbf{m} \not\leq \mathbf{N}}\rho_{\mathbf{m}}^2
%
\\
\nonumber 
& =  &
\big(  2 d^2 \eta_1^{2d-2}\eta^2_2 + 2\di\frac{ \eta_1^{2d}}{N_{max}^4}\big)\di \frac{ N_{max}^{2d+4}}{n} +
\di \sum_{\mathbf{m} \not\leq \mathbf{N}}\rho_{\mathbf{m}}^2\\
\label{inegal100}
&\leq & \di \frac{\eta N_{max}^{2d+4}}{n} + \di \sum_{\mathbf{m} \not\leq \mathbf{N}}\rho_{\mathbf{m}}^{2}\left(1+\sum_{i=1}^{d}m_{i}^{\beta_i}\right)\left(1+\sum_{i=1}^{d}m_{i}^{\beta_i}\right)^{-1},
\e*
with 
$\eta =  2d^2\eta_1^{2d-2}\eta_2^2 + {2\eta_1^{2d}}$.  





%



Observe that we can decompose 
$$\left\{\bm m\in \N^{d} \,/\, \bm m\not\leq \bm N\right\}= \bigcup\limits_{j=1}^{d}\left\{ \bm m\in \N^{d} / m_{j}>N_j {\rm \ and \ } \m m_{j+} \leq  \bm N_{j+}\right\},$$
where $\m m_{j+}=(m_{j+1},\cdots, m_d)$ and $\m N_{j+}=(N_{j+1},\cdots, N_d)$, with the convention 
$\m m_{d+} =  \m N_{d+} = \varnothing$. We also write 
$\m m_{j-}=(m_{1},\cdots, m_{j-1})$ and $\m N_{j-}=(N_{1},\cdots, N_{j-1})$, with the convention 
$\m m_{1-} =  \m N_{1-} = \varnothing$. 
%
We use this decomposition as follows:  
\be*
\di \E (\| \w{c}^{[{\mathbf{N}}]}-c\|_2^2 )
&\leq &
\di \frac{\eta N_{max}^{2d+4}}{n} + \di \sum_{j=1}^{d}\sum_{m_j>N_j}
\sum_{\begin{array}{l}
\tiny \m m_{j-} \in \N^{j-1}, 
\\
\tiny \m m_{j+} \leq  \bm N_{j+}
\end{array}}
\rho_{\mathbf{m}}^{2}\left(1+\sum_{i=1}^{d}m_{i}^{\beta_i}\right)\left(1+\sum_{i=1}^{d}m_{i}^{\beta_i}\right)^{-1}\\
&\leq &
\di \frac{\eta N_{max}^{2d+4}}{n} + \di \sum_{j=1}^{d}\sum_{m_j>N_j}
\sum_{\begin{array}{l}
\tiny \m m_{j-} \in \N^{j-1}, 
\\
\tiny \m m_{j+} \leq  \bm N_{j+}
\end{array}}
\rho_{\mathbf{m}}^{2}\left(1+\sum_{i=1}^{d}m_{i}^{\beta_i}\right)m_j^{-\beta_j}, 
\e*
 since  for all  $j=1,\cdots,d,$
 \be*
 \left(1+\sum_{i=1}^{d}m_{i}^{\beta_i}\right)^{-1}&\leq& m_{j}^{-\beta_j}.
\e*
Observe also that
\be*
b=\left(\sum_{j=1}^{d}\beta_{j}^{-1}\right)^{-1}\leq \beta_j, \text{ for all } j=1,\cdots,d.
\e*
Hence
\beg
\di \E (\| \w{c}^{[{\mathbf{N}}]}-c\|_2^2 )
&\leq &
\di \frac{\eta N_{max}^{2d+4}}{n} + \di \sum_{j=1}^{d}
\sum_{m_j>N_j} N_j^{-\beta_j}
\sum_{\begin{array}{l}
\tiny  \m m_{j-} \in \N^{j-1}, 
\nonumber
\\
\tiny  \m m_{j+} \leq  \bm N_{j+}
\end{array}}
\rho_{\mathbf{m}}^{2}\left(1+\sum_{i=1}^{d}m_{i}^{\beta_i}\right)
\nonumber
\\
&\leq &
\di \frac{\eta N_{max}^{2d+4}}{n} + \di d N_{max}^{-b} \sum_{\m N \in \N^d}
\rho_{\mathbf{m}}^{2}\left(1+\sum_{i=1}^{d}m_{i}^{\beta_i}\right)
\nonumber
\\
& \leq & 
\label{eq:rate} 
\di \frac{\eta N_{max}^{2d+4}}{n} +
L d N_{max}^{-b}.
\en
%
%
%
Since $N_{max}=[n^{\frac{1}{2d+b+4}}]$, we get
\be*
\E (\| \w{c}^{[{\mathbf{N}}]}-c\|_2^2 )
&\leq &
\di (\eta+Ld)n^{\frac{-b}{b+2d+4}}
%
\e*
which gives the result. 

%

%
%
\subsection*{Proof of Proposition  \ref{coellip}.}



According to (\ref{eq:rate}) we have
\be*
\di \E (\| \w{c}^{[{\mathbf{N}}]}-c\|_2^2 )
&\leq &
\di \frac{\eta N_{max}^{2d+4}}{n} +
L d N_{max}^{-b},
\e*
and it follows that 
\be*
\sup_{\di c\in \mathlarger{\mathcal{F}}_{\boldsymbol{\beta}}(L)} \di \E (\| \w{c}^{[{\mathbf{N}}]}-c\|_2^2 )
&\leq &
\di \frac{\eta N_{max}^{2d+4}}{n} +
L d N_{max}^{-b}. 
\e*
Since this inequality is satisfied for all $\m N \neq \m 0$ we obtain 
\be*
 \di \mathcal{R}^{\star}( \mathcal{F}_{\boldsymbol{\beta}}(L))
&\leq &
\di \frac{\eta N_{max}^{2d+4}}{n} +
L d N_{max}^{-b}, 
\e*
where the left hand side does not depend of $\m N$. 
By minimizing the right hand side, we obtain a minimum for 
$ \di N_{max}= n^{\frac{1}{2d+b+4}}\left(\frac{Lbd}{\eta (2d+4)}\right)^{\frac{1}{2d+b+4}}$ and  we conclude that
\be*
 \di \mathcal{R}^{\star}( \mathcal{F}_{\boldsymbol{\beta}}(L))\leq K v_n, \text{ where } K=\frac{Ld(2d+b+4)}{2d+4}\left(\frac{Lbd}{\eta(2d+4)}\right)^\frac{-b}{2d+b+4}.
\e*










%

\subsection{Proof of Proposition \ref{prop 5}}
We have
\be*
&\di \w{LSCV}({\mathbf N}) =& \di \int_{I^{d}} \left(\w{c}^{[{\mathbf N}]}({\mathbf u})\right)^2 d{\mathbf u} -\di \frac{2}{n}\di \sum_{i=1}^{m}\w{c}^{[{\mathbf N}]}_{-i}\left(\w F_{1}(X_{i1}),\cdots,\w F_{d}(X_{id})\right)\\
&&\hspace*{-0.7cm}= \sum_{{\mathbf m}\leq {\mathbf N}}\w{\rho}_{{\mathbf m}}^2 -\di \frac{2}{n}\sum_{i=1}^{n}\sum_{{\mathbf m}\leq {\mathbf N}}\w{\rho}_{{\mathbf m}}^{(-i)}\prod_{j=1}^{d}Q_{m_j}\Big(\w F_{j}(X_{ij})\Big)\\
&& \hspace*{-0.7cm}=
\sum_{{\mathbf m}\leq {\mathbf N}}\left(\frac{1}{n}\di\sum_{i=1}^{n}\prod_{j=1}^{d}Q_{m_j}\Big(\w F_{j}(X_{ij})\Big)\right)^2
\\
&& -\frac{2}{n}\sum_{i=1}^{n}\left( \sum_{{\mathbf m}\leq {\mathbf N}}\left(\frac{1}{n-1}\sum_{\substack{k=1,\\k\not\neq i }}^{n
}\prod_{j=1}^{d}Q_{m_j}\Big( \w F_{j}(X_{kj})\Big)\right)\prod_{j=1}^{d}Q_{m_j}\Big(\w F_{j}(X_{ij})\Big)\right)\\
&&\hspace*{-0.7cm}= \frac{1}{n^2}\sum_{{\mathbf m}\leq N}\left(\di\sum_{i=1}^{n}\prod_{j=1}^{d}Q_{m_j}^{2}\Big( \w F_{j}(X_{ij})\Big)+
\sum_{\substack{k\not\neq i }}\prod_{j=1}^{d}Q_{m_j}\Big(\w F_{j}(X_{kj})\Big)Q_{m_j}\Big(\w F_{j}(X_{ij})\Big)\right) \\
&&-\frac{2}{n(n-1)} \sum_{{\mathbf m}\leq {\mathbf N}}\sum_{\substack{k\not\neq i }}\prod_{j=1}^{d}Q_{m_j}\Big(\w F_{j}(X_{kj})\Big)Q_{m_j}\Big(\w F_{j}(X_{ij})\Big)\\
&&\hspace*{-0.7cm}= \di \frac{1}{n^2}\sum_{{\mathbf m}\leq N}\left(\di\sum_{i=1}^{n}\prod_{j=1}^{d}Q_{m_j}^{2}\Big(\w F_{j}(X_{ij})\Big)-\di \frac{n+1}{n-1}\sum_{k\neq i}\prod_{j=1}^{d}Q_{m_j}\Big(\w F_{j}(X_{ij})\Big)Q_{m_j}\Big(\w F_{j}(X_{kj})\Big) \right)
\e*
\subsection{Proof of Proposition \ref{prop 6}}
If the margins are known, we have
\be*
\di && \E\left(LSCV(N)\right) =  \E \left(\di \int_{I^{d}} \left(\tilde{c}^{[{\mathbf N}]}({\mathbf u})\right)^2 \mu(d{\mathbf u})\right) -\di 2\E\left(\di \frac{1}{n} \sum_{i=1}^{n}\tilde{c}^{[{\mathbf N}]}_{-i}\left(F_{1}(X_{i1}),\cdots,F_{d}(X_{id})\right)\right).
\e*
In addition we have
\be*
\di \E\left(\frac{1}{n}\sum_{i=1}^{n}\tilde{c}^{[{\mathbf N}]}_{-i}\left(F_{1}(X_{i1}),\cdots,F_{d}(X_{id})\right)\right)=
\di \frac{1}{n}\sum_{i=1}^{n}\sum_{m\leq N}\E\left(\tilde{\rho}_{m}\right)\E\left(\prod_{j=1}^{d}Q_{m_j}(F_{j}(X_{ij}))\right)
\e*
\be*
&=& \sum_{m\leq N}\E(\tilde{\rho}_{m})\int_{\R^{d}}\prod_{j=1}^{d}Q_{m_j}(F_{j}(x_{ij}))f(x_1,\cdots,x_d)\mu(\mathrm{d}{\mathbf x})\\
&=& \sum_{n\leq N}\E(\tilde{\rho}_n)\int_{I^{d}}\prod_{j=1}^{d}Q_{m_j}(u_{j})c\left(u_{1},\cdots ,u_{d})\right)\mu(\mathrm{d}{\mathbf u})\\
&=& \E \int_{I^{d}}\sum_{m\leq N}\tilde{\rho}_{m}\prod_{j=1}^{d}Q_{m_j}(u_j)c(u_1,\cdots,u_d)\mu(\mathrm{d}{\mathbf u}) \\
&=& \E\int_{I^{d}}\tilde{c}^{[{\mathbf N}]}({\mathbf u})c({\mathbf u})\mu(\mathrm{d}{\mathbf u})
\e*
which implies that
\be*
 \E\left(LSCV(N)\right) &= & \E \left(\di \int_{I^{d}} \left(\tilde{c}^{[{\mathbf N}]}({\mathbf u})\right)^{2}\mu(\mathrm{d}{\mathbf u}) -
2\int_{I^{d}}\tilde{c}^{[{\mathbf N}]}({\mathbf u})c({\mathbf u})\mu(\mathrm{d}{\mathbf u})\right)\\
&=&\E\left(\int_{I^{d}}\left(\tilde{c}^{[{\mathbf N}]}({\mathbf u})-c({\mathbf u})\right)^{2}\mu(\mathrm{d}{\mathbf u}) \right) - \int_{I^{d}}(c({\mathbf u}))^{2}\mu(\mathrm{d}{\mathbf u}),
\e*
and we obtain
\beg\label{lscv_know_margins}
\E\left(LSCV(N)\right)&=& \E\|\tilde{c}^{[{\mathbf N}]}-c\|_{2}^{2} -\|c\|_{2}^{2}.
\en
If the margins are unknown, we combine Corollary \ref{relation_mise} and equation (\ref{lscv_know_margins}) to obtain the result.

\section{Modifying the copula density to be square integrable}
\label{appendL2}
If assumption (\ref{condition}) is not satisfied, one can modify the copula density $c$ with a shrinkage  function, 
writing 
\be*
c(\boldsymbol{\theta},\m u)  & = &
s(\boldsymbol{\theta}, \m u)   c(\m u), 
\e*
for some $\boldsymbol{\theta} =(\theta_1,\cdots, \theta_p)$. 
For instance, $s$ can be an exponential tilting defined as follows 
$$
\di s(\boldsymbol{\theta}, \m u)=\exp\left(- \di\sum_{j=1}^p \di\frac{\theta_j}{ u_j} \right) := \di\prod_{j=1}^d s_j(\theta_j,u_j), $$

 {\rm with \ } $\theta_j > 0$. 
In that case we obtain
\begin{align*}
c(\boldsymbol{\theta}, \m u)
& = \di\sum_{\m m \in \N^p} \rho_{\boldsymbol{\theta},\m m}\mathbf{Q}_{\mathbf{m}}({\mathbf u})),
\end{align*}
where
\begin{align*}
    \rho_{\boldsymbol{\theta},\m m} &=
    \E\left(\prod_{j=1}^{d}s_j(\theta_j,U_j)Q_{m_j}(U_j)\right).
\end{align*}
We then modify our estimators as follows
\begin{align}
\label{modify}
\w \rho_{\boldsymbol{\theta}, \m j} & =  \di\frac{1}{n} \di\sum_{i=1}^{n} \left(\prod_{j=1}^{d}s_j(\theta_j,\w U_{i,j})Q_{m_j}(\w U_{i,j})\right) 
\end{align} 
and we get a $\mathbf N$-th order  estimator of $c(\boldsymbol{\theta},\m u)$ as 
\be*
\di \w{c}^{[{\mathbf N}]}(\boldsymbol{\theta},\m u) &=& \sum_{{\mathbf m}\leq {\mathbf N}}\w{\rho}_{\boldsymbol{\theta}, {\mathbf m}}\prod_{j=1}^{d}Q_{m_j}(u_j).
\e*
Clearly the positive functions $c(\boldsymbol{\theta},.)$ do not satisfy the properties of copula densities but we expect  that the transformation made them square integrable. 
It gives a way to rewrite the main results of the paper. Since we apply a known transformation  we can 
find an estimation of   $c(.)$ thanks to that   of 
$c(\boldsymbol{\theta},.)$.
More precisely we use 
\be*
\w c^{[\m N]}(\m u ) 
& = & \di s^{-1}(\m \theta, \m u) \w c^{[\m N]}(\boldsymbol{\theta},\m u).
\e*
The fact that all factors $\di s_j(\m \theta_j,U_j)$ and $\di s_j(\m \theta_j,\w U_j)$ are bounded is essential here to generalize our methodology and to get the same asymptotic results.

We illustrate this shrinking approach through the bivariate Clayton and Gumbel copulas as follows.   
\paragraph{The bivariate Clayton copula.} 
Its copula density is given by
\be* 
c(u_1,u_2) 
& = &
(\beta+1) (u_1u_2)^{-(\beta+1)} \big(u_1^{-\beta} + u_2^{-\beta} -1 \big)^{-(2\beta+1)/\beta},
\e*
with $u_1,u_2 \in (0,1)$ and for  $\beta > 0$.  
We choose 
$$
s(\boldsymbol{\theta} , \m u)=\exp\left(-\di\frac{\theta_1}{ u_1} -\di\frac{\theta_2}{  u_2}\right),
$$
and an easy computation shows that 
\be* 
\di\int_0^1 \di\int_0^1 \left( 
s(\boldsymbol{\theta} , \m u)c(u_1,u_2)\right)^2 du_1du_2
& \leq  & (\beta+1)
\di\int_0^1 \di\int_0^1 \exp\left(-\di\frac{2\theta_1}{ u_1} -\di\frac{2\theta_2}{  u_2}\right) (u_1u_2)^{-2(\beta+1)} du_1du_2
\\
 & = & 
\di\int_1^{\infty} \di\int_1^{\infty} \exp\left(-{2{\theta_1} x} - {2\theta_2 y}\right) (xy)^{2\beta} dxdy
\\ 
& < \infty& {\rm \ \ for \ all \ } \theta_1>0,\theta_2 >0.
\e*

\paragraph{The bivariate Gumbel copula.}  Its density is given by 
\be*
c(u_1,u_2) 
& = & 
\di\frac{C(u_1,u_2)}{u_1u_2} 
\left( 
\ell_{\beta}(u_1) + \ell_{\beta}(u_2)\right)^{-2+2/\beta} 
\ell_{\beta-1}(u_1) \ell_{\beta-1}(u_2)
\left(
1+(\beta-1)\big(\ell_{\beta}(u_1)+\ell_{\beta}(u_2)\big)^{-1/\beta}
\right) 
\e*
where 
\be*
C(u_1,u_2)  = 
\exp\left(-\big(\ell_{\beta}(u_1)+\ell_{\beta}(u_2)\big)^{1/\beta}\right)
& {\rm and }  & 
\ell_{\beta}(u) = \big(-\log(u)\big)^{\beta},
\e*
with $u_1,u_2 \in (0,1)$ and $\beta \geq 1$. 
We can chek that 
\be*
\di\int_0^1 \di\int_0^1 \left( u_{1}^{\theta_1}u_{2}^{\theta_2} c(u_1,u_2)\right)^2 du_1du_2 & < & \infty {\rm \ \ for \ all \ } \theta_1>0,\theta_2 >0. 
\e*
In that case the shrinkage factor can be simply $s(\boldsymbol{\theta}, \m u) = u_{1}^{\theta_1} u_{2}^{\theta_2}$.

In conclusion, since $s(\boldsymbol{\theta},\m u)$ is known we can choose arbitrary $\boldsymbol{\theta}$ such that for a given small $\epsilon >0$ we have 
\be* 
1-\epsilon \leq \di\frac{c(\boldsymbol{\theta},\m u)}{c(\m u)}:=  s(\boldsymbol{\theta},\m u) < 1.   
\e*
The proximity between $c(.)$ and $c(\boldsymbol{\theta},)$ explains the very good behavior of the method even for the case where (\ref{condition}) is not satisfied. In that case 
the coefficient estimators 
$\w \rho_{\m j}$ and 
$\w \rho_{\boldsymbol{\theta}, \m j}$ are very close and so are the density estimators $\w c^{[\m N]}(.)$ and $\w c^{[\m N]}(\boldsymbol{\theta},)$

We illustrate numerically this remark 
in the six sample case  for  the Clayton and Gumbel copulas, using a  tilting exponential factor $s(\boldsymbol{\theta} , \m u)=\exp\left(-\di\frac{\theta_1}{ u_1} -\di\frac{\theta_2}{  u_2}\right)$, with $(\theta_1,\theta_2)=(0.001,0.001)$ arbitrary small such that  for all $\m u \in I^2$: 
$0.998 < s(\boldsymbol{\theta} , \m u)< 1$.  
The results are given in Table \ref{modifErros}. As expected we can observe that the accuracies of these new estimators are very similar to the accuracies of those without shrinkage. 
This means that our estimator gives a very good estimate of the function $c(\boldsymbol{\theta},.)$ that is itself very close to the desired copula density $c$. 
And we could still improve  these  estimations by taking a smaller value for $\theta$.



\newpage
\section*{TABLES AND FIGURES}

\begin{table}[!h]
\centering
\resizebox{18cm}{8cm}{%
\begin{tabular}{|c|c|c||c|c|c|c|c|c|c|}
\hline
\multirow{21}{*}{\rotatebox{90}{sample size $n=500$ }}
&\multicolumn{2}{|c}{\multirow{2}{*}{Copulas} }& \multicolumn{7}{|c|}{Methods} \\ \cline{4-10}
&\multicolumn{2}{|c|}{} & Emp & Beta & Chek & Berns10  & Berns25 & CN  & $N_{opt}$ \\\hhline{~=========} 
&\multirow{7}{*}{\rotatebox{90}{$\tau=0.3$ }}
& Gauss & 1.78(\underline{0.44}) & 1.63(0.47) &    1.75(0.44) & 2.62(0.93) & $\bm{1.51}$(0.68) & 1.53(0.50) & 2  \\  
& & Frank & 1.76(\underline{0.41}) & 1.63(0.44) & 1.74(0.42) & 2.58(0.87) & 1.50(0.67) & $\bm{1.29}$(0.45) & 1 \\   
& & Student17 & 1.80(\underline{0.45}) & 1.66(0.48) & 1.78(0.45) & 2.62(0.95) & $\bm{1.53}$(0.70) & 1.55(0.52) & 2  \\
& & Gumbel & 1.81(\underline{0.45}) & 1.68(0.49) & 1.79(0.45) & 2.60(0.94) & 1.56(0.68) & $\bm{1.50}$(0.54) & 3 \\
& & Joe & 1.80(\underline{0.47}) & 1.67(0.50) & 1.79(0.48) & 2.74(0.92) & 1.65(0.73) & $\bm{1.61}$(0.53) & 7 \\
& & Clayton & 1.79(0.47) & 1.64(0.49) & 1.76(0.47) & 2.68(0.88) & 1.60(0.68) & $\bm{1.53}$(\underline{0.44}) & 5 \\
\hhline{~=========}
& \multirow{7}{*}{\rotatebox{90}{$\tau=0.55$ }}
& Gauss & 1.14(\underline{0.23}) & 0.98(0.26) & 1.07(0.24) & 4.02(0.52) & 1.67(0.55) & $\bm{0.93}$(0.29) & 5 \\  
& & Frank & 1.15(\underline{0.23}) & 1.02(0.25) & 1.11(0.23) & 3.97(0.47) & 1.69(0.48) & $\bm{0.89}$(0.29) & 3  \\  
& & Student17 & 1.17(\underline{0.25}) & 1.01(0.27) & 1.11(\underline{0.25})& 4.0(0.54) & 1.68(0.57) & $\bm{0.96}$(0.30) & 5  \\
& & Gumbel & 1.18(0.25) & 1.04(0.28) & 1.13(\underline{0.25}) & 3.97(0.57) & 1.68(0.30) & $\bm{0.99}$(0.30) & 6 \\
& & Joe & 1.58(0.26) & $\bm{1.02}$(0.29) & 1.11(\underline{0.26}) & 3.96(0.57) & 1.75(0.56) & $\bm{1.03}$(0.29) & 10 \\
& & Clayton & 1.17(\underline{0.27}) & $\bm{1.01}$(0.29) & 1.09(\underline{0.27}) & 3.93(0.57) & 1.72(0.56) & $\bm{1.03}$(0.30) & 9  \\
\hhline{~=========}
&\multirow{7}{*}{\rotatebox{90}{$\tau=0.8$ }}
& Gauss & 0.31(\underline{0.044}) & $\bm{0.23}$(0.051) & 0.26(0.044) & 4.76(0.092) & 1.97(0.098) & 0.27(0.050) & 10 \\  
& & Frank & 0.51(\underline{0.058}) &$\bm{0.33}$(0.072) & 0.38(0.057) & 4.75(0.10) & 1.98(0.11) & 0.41(0.073) & 8  \\  
& & Student17 & 0.52(\underline{0.064}) & $\bm{0.33}$(0.083) & 0.38(0.064) & 4.76(0.14) & 1.98(0.15) & 0.41(0.079) & 10 \\
& & Gumbel & 0.54(0.069) & 0.36(0.085) & $\bm{0.41}$(\underline{0.067}) & 4.75(0.15) & 1.97(0.17) & 0.43(0.083) & 10  \\
& & Joe & 0.52(0.081) & $\bm{0.36}$(0.096) & 0.40(\underline{0.076}) & 4.71(0.17) & 1.97(0.18) & 0.44(0.090) & 16  \\
& & Clayton & 0.54(\underline{0.077}) & $\bm{0.35}$(0.099) & 0.39(0.079) & 4.70(0.17) & 1.96(0.18) & 0.46(0.085) & 17 \\ \hhline{~=========}
& \multicolumn{2}{|c||}{Independent} & 2.30(0.59) & 2.14(0.63) & 2.28(0.59) & 1.30(0.62) & 1.64(0.66) & $\bm{4.81\text{e-}32}$(\underline{0.00}) & 0 \\
\hhline{==========}
\multirow{21}{*}{\rotatebox{90}{sample size $n=1000$ }}
&\multicolumn{2}{|c}{\multirow{2}{*}{Copulas} }& \multicolumn{7}{|c|}{Methods} \\ \cline{4-10}
&\multicolumn{2}{|c|}{} & Emp & Beta & Chek & Berns10 & Berns25 & CN & $N_{opt}$ \\\hhline{~=========} 
&\multirow{7}{*}{\rotatebox{90}{$\tau=0.3$ }}
& Gauss & 1.25(0.32) & 1.18(0.34) & 1.24(\text{0.33}) & 2.61(0.67) & 1.26(0.57) & $\bm{1.02}$(0.39) & 3  \\  
& & Frank & 1.25(\underline{0.30}) & 1.82(0.31) & 1.24(\underline{0.30}) & 2.57(0.63) & 1.28(0.54) & $\bm{0.87}$(0.39) & 2 \\  
& & Student17 & 1.26(\underline{0.33}) & $\bm{1.20}$(0.34) & 1.26(\underline{0.33}) & 2.60(0.68) & 1.27(0.57) & 1.26(0.34) & 2  \\
& & Gumbel & 1.28(\underline{0.32}) & 1.21(0.34) & 1.27(\underline{0.32}) & 2.59(0.70) & 1.30(0.56) & $\bm{1.11}$(0.37) & 5 \\
& & Joe & 1.23(\underline{0.30}) & 1.16(0.32) & 1.16(0.31) & 2.59(0.63) & 1.29(0.53) & $\bm{1.09}$(0.34) & 7 \\
& & Clayton & 1.24(\underline{0.31}) & 1.18(0.32) & 1.24(\underline{0.31}) & 2.58(0.65) & 1.30(0.51) & $\bm{1.16}$(0.33) & 12 \\
\hhline{~=========}
&\multirow{7}{*}{\rotatebox{90}{$\tau=0.55$ }}
& Gauss & 79.0(\underline{17.3}) & 72.0(18.0) & 77.0(17.3) & 402.0(37.0) & 165.0(40.0) & $\bm{67.0}$(20.0) & 6 \\  
& & Frank & 0.81(0.16) & 0.75(0.17) & 0.79(\underline{0.16}) & 3.97(0.34) & 1.67(0.36) & $\bm{0.67}$(0.19) & 3  \\  
& & Student17 & 0.80(\underline{0.18}) & 0.74(0.19) & 0.79(0.18)  & 4.00(0.38) & 1.60(0.42) & $\bm{0.69}$(0.21) & 5 \\
& & Gumbel & 0.82(0.17) & 0.76(0.18) & 0.81(\underline{0.17}) & 3.98(0.40) & 1.66(0.43) &  $\bm{0.75}(0.19)$ & 12 \\
& & Joe & 0.79(0.18) & $\bm{0.73}$(0.18) & 0.77(\underline{0.174}) & 3.91(0.39) & 1.64(0.41) & 0.74(0.18) & 17 \\
& & Clayton & 0.79(0.18) & $\bm{0.72}$(0.19) & 0.77(\underline{0.18}) & 3.91(0.39) & 1.64(0.41) & 0.73 (0.19) & 13  \\
\hhline{~=========}
&\multirow{7}{*}{\rotatebox{90}{$\tau=0.8$ }}
& Gauss & 0.31(0.04) & $\bm{0.23}(0.05)$ & 0.26(\underline{0.04}) & 4.76(0.09) & 1.97(0.10) & 0.27(0.05) & 14 \\  
& & Frank & 0.32(\underline{0.04}) & $\bm{0.25}$(0.05) & 0.27(\underline{0.04}) & 4.74(0.07) & 1.97(0.08) & 0.28(0.05) & 9  \\  
& & Student17 & 0.32(\underline{0.04}) & $\bm{0.24}$(0.05) & 0.27(0.05) & 4.75(0.10) & 1.97(0.10) & 0.28(0.05) & 13  \\
& & Gumbel & 0.34(\underline{0.05}) & $\bm{0.26}$(0.06) & 0.29(\underline{0.05}) & 4.75(0.11) & 1.97(0.12) & 0.30(0.06) & 16  \\
& & Joe & 0.33(0.05) & 0.25(0.06) & $\bm{0.28}$(\underline{0.05}) & 4.70(0.14) & 1.95(0.13) & 0.29(0.06) & 20  \\
& & Clayton & 0.33(0.05) & $\bm{0.25}$(0.06) & 0.28(\underline{0.05}) & 4.70(0.11) & 1.96(0.12) & 0.30(0.06) & 20 \\
\hhline{~=========}
& \multicolumn{2}{|c||}{Independent } & 1.62(0.43) & 1.55(0.45) & 1.62(0.44) & 0.92(0.45) & 1.16(0.49) & $\bm{4.81\text{e-}15}$(\underline{0.00}) & 0 \\
\hline
\end{tabular}
}
\caption{\label{MIAE-copula} Bivariate copulas: Relative MIAE $\times 100$.
Bold values show the minimun of the relative MIAE and the
underline the minimun of standard deviation for the corresponding
copula.}
\end{table}

\FloatBarrier

\begin{table}[!h]
\centering
\resizebox{18cm}{8cm}{%
\begin{tabular}{|c|c|c||c|c|c|c|c|cc|}
\hline
\multirow{21}{*}{\rotatebox{90}{sample size $n=500$ }}
&\multicolumn{2}{|c}{\multirow{2}{*}{Copulas} }& \multicolumn{7}{|c|}{Methods} \\ \cline{4-10}
&\multicolumn{2}{|c|}{} & Emp & Beta & Chek &  Berns10 & Berns25 & CN  & $N_{opt}$ \\\hhline{~=========} 
&\multirow{7}{*}{\rotatebox{90}{$\tau=0.3$ }}
& Gauss & 3.76(2.15) & 3.26(2.14) & 3.71(2.17) & 2.78(4.33) & 2.78(2.55) & $\bm{2.58}$(\underline{2.01}) & 2 \\  
& & Frank & 3.61(1.88) & 3.13(1.88) & 3.57(1.89) & 7.17(4.46) & 7.17(2.63) & $\bm{1.77}$(\underline{0.16}) & 1 \\   
& & Student17 & 3.82(2.18) & 3.32(2.18) & 3.77(2.20) & 6.62(4.48) & 2.85(2.66) & $\bm{2.66}$(\underline{2.08}) & 2  \\
& & Gumbel & 3.80(2.12) & 3.32(2.10) & 3.76(2.12) & 6.78(4.32) & 2.90(2.56) & $\bm{2.64}$(\underline{2.09})  & 3 \\
& & Joe & 3.77(\underline{2.20}) & 3.30(\underline{2.20}) & 3.72(2.21) & 8.06(4.64) & 3.22(2.84) & $\bm{3.06}$(2.21) & 7 \\
& & Clayton & 3.89(2.29) & 3.38(\underline{2.27}) & 3.82(2.29) & 7.92(4.62) & 3.17(2.28) & $\bm{2.92}$(2.28) & 5 \\
\hhline{~=========}
&\multirow{7}{*}{\rotatebox{90}{$\tau=0.55$ }}
& Gauss & 1.86(0.92) & 1.53(0.92) & 1.81(0.92) & 16.0(4.22) &3.45(2.16) & $\bm{1.23}$(\underline{0.90}) & 5 \\  
& & Frank & 1.79(0.82) & 1.48(0.83) & 1.75(0.81) & 18.2(3.99) & 3.99(2.12) & $\bm{0.99}$(\underline{0.73}) & 3  \\  
& & Student17 & 1.91(0.96) & 1.58(0.97) & 1.86(0.96) & 16.0(4.36) & 3.62(2.26) & $\bm{1.28}$(\underline{0.95}) & 5  \\
& & Gumbel & 1.91(0.94) & 1.59(0.94) & 1.86(\underline{0.93})  & 17.0(4.42) & 3.72(2.28) & $\bm{1.37}$(\underline{0.93}) & 6 \\
& & Joe & 1.88(1.01) & 1.58(1.01) & 1.84(\underline{1.00}) & 18.0(4.45) & 4.24(2.37) & $\bm{1.53}$(1.01) & 10 \\
& & Clayton & 1.99(1.14) & 1.66(1.15) & 1.92(1.14) & 4.32(4.72) & 4.32(2.52) & $\bm{1.59}$(\underline{1.13}) & 9  \\
\hhline{~=========}
&\multirow{7}{*}{\rotatebox{90}{$\tau=0.8$ }}
& Gauss & 0.21(0.84) & 0.17(0.88) & 0.20(0.08) & 0.32(1.15) & 6.66(0.64) & $\bm{0.15}$(\underline{0.08}) & 10 \\  
& & Frank & 0.43(0.14) & 0.31(0.15) & 0.38(\underline{0.13}) & 0.35(1.31) & 7.40(0.72) & $\bm{0.26}$(\underline{0.13}) & 8  \\  
& & Student17 & 0.47(0.17) & 0.33(0.19) & 0.42(0.17) & 0.33(1.72) & 6.78(0.96) & $\bm{0.29}$(\underline{0.16}) & 10\\
& & Gumbel & 0.49(0.17) & 0.35(0.19) & 0.44(0.17) & 0.33(1.84) & 6.87(1.02) & $\bm{0.32}$(\underline{0.16}) & 10  \\
& & Joe & 0.50(\underline{0.21}) & $\bm{0.38}$(0.23) & 0.45(0.20) & 0.34(1.88) & 7.50(1.05) & 0.38(\underline{0.21}) & 16  \\
& & Clayton & 0.51(\underline{0.21}) & $\bm{0.38}$(0.24) & 0.46(0.22) & 0.34(1.95) & 7.58(1.08) & 0.40(\underline{0.21}) & 17 \\
\hhline{~=========}
& \multicolumn{2}{|c||}{Independent}& 5.38(3.10) & 4.72(3.10) & 5.32(3.10) & 1.88(1.90) & 2.89(2.50) & $\bm{8.02\text{e-}29}$(\underline{0.00}) & 0 \\
\hhline{==========}
\multirow{21}{*}{\rotatebox{90}{sample size $n=1000$ }}
& \multicolumn{2}{|c}{\multirow{2}{*}{Copulas} }& \multicolumn{7}{|c|}{Methods} \\ \cline{4-10}
& \multicolumn{2}{|c|}{} & Emp & Beta & Chek & Berns10 & Berns25 & CN & $N_{opt}$\\\hhline{~=========} 
&\multirow{7}{*}{\rotatebox{90}{$\tau=0.3$ }}
& Gauss & 1.88(1.14) & 1.71(1.14) & 1.87(1.15) & 5.90(3.00) & 1.85(1.65) & $\bm{1.13}$(\underline{1.13}) & 3  \\  
& & Frank & 1.66(1.00) & 1.66(0.99) & 1.82(1.00) & 6.69(3.08) & 2.03(1.69) &  $\bm{0.94}$(\underline{0.89}) & 2 \\  
& & Student17 & 1.90(\underline{1.15}) & 1.73(\underline{1.15}) & 1.89(\underline{1.15}) & 5.99(\underline{1.15}) & 1.89(\underline{1.15}) & $\bm{1.67}$(\underline{1.15}) & 2  \\
& & Gumbel & 1.90(\underline{1.11}) & 1.74(\underline{1.11}) & 1.90(\underline{1.11}) & 6.23(3.06) & 1.94(1.66) & $\bm{1.46}$(\underline{1.11}) & 5 \\
& & Joe & 1.74(0.98) & 1.58(0.97) & 1.73(0.98) & 6.95(2.83) & 1.93(1.54) & $\bm{1.40}$(\underline{0.97}) & 7 \\
& & Clayton & 1.91(1.06) & 1.74(\underline{1.05}) & 1.90(1.06) & 7.03(3.0) & 2.03(1.56) & $\bm{1.67}$(\underline{1.05}) & 12 \\
\hhline{~=========}
& \multirow{7}{*}{\rotatebox{90}{$\tau=0.55$ }}
& Gauss & 0.92(0.50) & 0.81(0.50) & 0.91(0.50) & 1.62(3.00) & 3.14(1.49) & $\bm{0.66}$(\underline{0.49}) & 6 \\  
& & Frank & 0.91(0.42) & 0.81(0.42) & 0.91(2.85) & 0.18(2.85) & 3.69(1.46) & $\bm{0.56}$(\underline{0.37}) & 3  \\  
& & Student17 & 0.94(0.50) & 0.83(0.50) & 0.93(4.99)  & 0.16(3.05) & 3.19(1.53) & $\bm{0.65}$(\underline{0.49}) & 5 \\
& & Gumbel & 0.96(\underline{0.47}) & 0.85(\underline{0.47}) & 0.95(\underline{0.47}) & 0.17(3.13) & 3.35(1.57) & $\bm{0.81}$(\underline{0.47}) & 12 \\
& & Joe & 0.90(0.46) & $\bm{0.80}$(0.46) & 0.89(0.46) & 0.18(2.94) & 3.63(1.51) & $\bm{0.80}$(\underline{0.46})  & 17 \\
& & Clayton & 0.02(\underline{0.49}) & 0.01(\underline{0.49}) & 0.95(\underline{0.49}) & 4.00(3.09) & 3.74(1.57) & $\bm{0.82}$(\underline{0.49}) & 13  \\
\hhline{~=========}
& \multirow{7}{*}{\rotatebox{90}{$\tau=0.8$ }}
& Gauss & 0.32(4.39) & 0.17(0.09) & 0.20(0.08) & 0.32(1.53) & 6.66(0.64) & $\bm{0.15}$(\underline{0.08}) & 14 \\  
& & Frank & 0.32(3.91) & 0.25(4.52) & 0.27(\underline{0.06}) & 0.34(0.90) & 7.33(0.50) & $\bm{0.13}$(\underline{0.06}) & 9  \\  
& & Student17 & 0.22(0.08) & 0.17(0.09) & 0.21(0.08) & 0.33(1.21) & 6.68(0.66) & $\bm{0.15}$(\underline{0.08}) & 13  \\
& & Gumbel & 0.23(\underline{0.09}) & 0.19(0.09) & 0.22(\underline{0.09}) & 0.33(1.29) & 6.82(0.71) & $\bm{0.18}$(\underline{0.09})& 16  \\
& & Joe & 0.23(0.09) &  $\bm{0.19}$(0.10) & 0.22(\underline{0.09}) & 0.34(1.26) & 7.32(0.71) & $\bm{0.19}$(0.09) & 20  \\
& & Clayton & 0.24(0.10) & 0.19(0.11) & 0.23(0.10) & 0.34(0.11) & 7.47(0.72) & $\bm{0.19}$(\underline{0.10})  & 20 \\
\hhline{~=========}
&\multicolumn{2}{|c||}{Independent }  & 2.71(43.0) & 2.48(1.6) & 2.70(2.0) & 0.96(29.0) & 1.46(34.0) & $\bm{8.02\text{e-}29}$(\underline{0.00}) & 0 \\
\hline
\end{tabular}
}
\caption{\label{MISE-copula} Bivariate copulas: Relative MISE $\times 10^4$.
Bold values show the minimun of relative MISE and the underline the minimun of standard deviation for the corresponding copula.}
\end{table}

\FloatBarrier

\begin{table}[!h]
\centering
\resizebox{18cm}{8cm}{%
\begin{tabular}{|c|c|c||c|c|c|c|c|cc|}
\hline
\multirow{21}{*}{\rotatebox{90}{sample size $n=500$ }}
&\multicolumn{2}{|c}{\multirow{2}{*}{Copulas} }& \multicolumn{7}{|c|}{Methods}  \\ \cline{4-10}
&\multicolumn{2}{|c|}{} & Emp & Beta & Chek &  Berns10 & Berns25 & CN  & $N_{opt}$ \\\hhline{~=========} 
&\multirow{7}{*}{\rotatebox{90}{$\tau=0.3$ }}
 & Gauss & 2.51(0.53) & 2.03(0.51) &  2.40(0.53) & 1.68(0.54) & 1.38(0.50) & $\bm{1.29}$(\underline{0.46}) & 2  \\  
& & Frank & 2.37(0.51) & 1.99(0.49) & 2.36(0.59) & 1.97(0.59) & 1.42(0.54) & $\bm{0.93}$(\underline{0.30}) & 1 \\   
& & Student17 & 2.41(0.53) & 2.03(0.51) & 2.39(0.53) & 1.72(0.53) & 1.39(0.51) & $\bm{1.30}$(\underline{0.46}) & 2  \\
& & Gumbel & 2.42(0.54) & 2.04(0.52) & 2.40(0.53) & 1.92(\underline{0.45}) & 1.42(0.50) & $\bm{1.41}$(0.46) & 3 \\
& & Joe & 2.39(0.53) & 2.02(0.52) & 2.37(0.53) & 2.42(\underline{0.40}) & $\bm{1.53}$(0.51) & 1.70(0.51) & 7 \\
& & Clayton & 2.41(0.47) & 1.64(0.49) & 17.6(0.47) & 2.68(0.88) & 1.60(0.68) & $\bm{1.53}$(\underline{0.54}) & 5 \\
\hhline{~=========}
&\multirow{7}{*}{\rotatebox{90}{$\tau=0.55$ }}
 & Gauss & 2.13(0.45) & 1.76(0.44) & 2.12(0.45) & 3.00(0.41) & 1.72(0.47) & $\bm{1.23}$(\underline{0.41}) & 5 \\  
& & Frank & 2.09(0.43) & 1.71(0.42) & 2.08(0.43) & 3.70(0.42) & 1.98(0.49) & $\bm{1.01}$(\underline{0.34}) & 3  \\  
& & Student17 & 2.13(0.45) & 1.75(0.42) & 2.11(0.45) & 3.05(0.42) & 1.75(0.48) & $\bm{1.23}$(\underline{0.41}) & 5  \\
& & Gumbel & 2.11(0.46) & 1.73(0.44) & 2.09(0.46) & 3.35(\underline{0.36}) & 1.86(0.44) & $\bm{1.32}$(0.43) & 6 \\
& & Joe & 2.10(0.46) & 1.74(0.45) & 2.09(0.46) & 4.11(\underline{0.28}) & 2.18(0.37) & $\bm{1.56}$(0.45) & 10 \\
& & Clayton & 2.15(0.48) & 1.77(0.46) & 2.13(0.48) & 4.05(\underline{0.30}) & 2.15(0.38) & $\bm{1.55}$(0.45) & 9  \\
\hhline{~=========}
& \multirow{7}{*}{\rotatebox{90}{$\tau=0.8$ }}
& Gauss & 1.09(0.22) & 0.91(0.22) & 1.08(0.22) & 5.30(\underline{0.12}) & 2.70(0.17) & $\bm{0.77}$(0.20) & 10 \\  
& & Frank & 1.51(0.30) & 1.21(0.30) & 1.49(0.30) & 6.28(\underline{0.13}) & 3.31(0.21) & $\bm{0.84}$(0.25) & 8  \\  
& & Student17 & 1.55(0.31) & 1.21(0.31) & 1.53(0.31) & 5.40(\underline{0.18}) & 2.81(0.24) & $\bm{0.94}$(0.28) & 10\\
& & Gumbel & 1.45(0.31) & 1.21(0.30) & 1.53(0.30) & 5.61(\underline{0.15}) & 2.96(0.22) & $\bm{0.96}$(0.28) & 10  \\
& & Joe & 1.50(0.33) & 1.25(0.32) & 1.50(0.33) & 6.40(\underline{0.11}) & 3.50(0.14) & $\bm{1.20}$(0.33) & 16  \\
& & Clayton & 1.60(0.33) & 1.28(0.34) & 1.58(0.34) & 6.50(\underline{0.11}) & 3.53(0.15) & $\bm{1.25}$(0.33) & 17 \\
\hhline{~=========}
& \multicolumn{2}{|c||}{Independent}& 2.49(0.56) & 2.11(0.53) & 2.48(0.56) & 0.96(0.40) & 1.31(0.47) & $\bm{1.17\text{e-}14}$(\underline{0.00}) & 0 \\
\hhline{==========}
\multirow{21}{*}{\rotatebox{90}{sample size $n=1000$ }}
&\multicolumn{2}{|c}{\multirow{2}{*}{Copulas} }& \multicolumn{7}{|c|}{Methods}  \\ \cline{4-10}
&\multicolumn{2}{|c|}{} & Emp & Beta & Chek &  Berns10 & Berns25 & CN & $N_{opt}$\\\hhline{~=========} 
&\multirow{7}{*}{\rotatebox{90}{$\tau=0.3$ }}
& Gauss & 1.70(0.38) & 1.51(0.37) & 1.70(0.38) & 1.55(0.39) & 1.08(0.38) & $\bm{0.99}$(\underline{0.35}) & 3  \\  
& & Frank & 1.68(0.37) & 1.49(0.36) & 1.68(0.37) & 1.90(0.43) & 1.18(0.42) & $\bm{0.75}$(\underline{0.32}) & 2 \\  
& & Student17 & 1.69(0.38) & 1.51(0.37) & 1.69(0.38) & 1.59(0.39) & 1.09(0.39) & $\bm{1.07}$(\underline{0.32}) & 2  \\
& & Gumbel & 1.70(0.37) & 1.51(0.36) & 1.70(0.37) & 1.82(\underline{0.31}) & 1.14(0.37) & $\bm{1.11}$(0.35) & 5 \\
& & Joe & 1.66(0.36) & 1.47(0.35) & 1.65(0.36) & 2.33(\underline{0.24}) & 1.24(0.31) & $\bm{1.18}$(0.34) & 7 \\
& & Clayton & 1.68(0.37) & 1.49(0.36) & 1.68(0.37) & 1.90(0.43) & 1.18(0.42) & $\bm{0.75}$(\underline{0.32}) & 12 \\
\hhline{~=========}
&\multirow{7}{*}{\rotatebox{90}{$\tau=0.55$ }}
& Gauss & 1.49(0.32) & 1.30(0.31) & 1.49(0.33) & 2.92(0.31) & 1.53(0.34) & $\bm{0.93}$(\underline{0.30}) & 6 \\  
& & Frank & 1.49(0.32) & 1.31(0.31) & 1.49(0.32) & 3.67(0.31) & 1.85(0.36) & $\bm{0.79}$(\underline{0.24}) & 3  \\  
& & Student17 & 1.49(0.33) & 1.29(0.32) & 1.48(0.33) & 2.95(0.31) & 1.56(0.34) & $\bm{0.87}$(\underline{0.29}) & 5\\
& & Gumbel & 1.50(0.31) & 1.31(0.31) & 1.49(0.32) & 3.29(\underline{0.26}) & 1.71(0.31) & $\bm{1.17}$(0.30) & 12 \\
& & Joe & 1.47(0.32) & 1.28(0.31) & 1.47(0.32) & 4.06(\underline{0.19}) & 2.03(0.23) & $\bm{1.23}$(0.31) & 17 \\
& & Clayton & 1.52(0.33) & 1.33(0.32) & 1.51(0.33) & 4.00(\underline{0.19}) & 2.00(0.24) & $\bm{1.21}$(0.32) & 13  \\
\hhline{~=========}
& \multirow{7}{*}{\rotatebox{90}{$\tau=0.8$ }}
& Gauss & 1.09(0.22) & 0.91(0.22) & 1.08(0.22) & 5.35(\underline{0.12}) & 2.69(0.17) & $\bm{0.77}$(0.20) & 14 \\  
& & Frank & 1.06(0.20) & 0.91(0.20) & 1.00(0.20) & 6.27(\underline{0.093}) & 3.27(0.15) & $\bm{0.64}$(0.18) & 9  \\  
& & Student17 & 1.09(0.22) & 0.90(0.21) & 1.08(0.22) & 5.38(\underline{0.13}) & 2.72(0.17) & $\bm{0.74}$(0.20) & 13  \\
& & Gumbel & 1.10(0.23) & 0.93(0.23) & 1.10(0.23) & 5.60(\underline{0.11}) & 2.90(0.16) & $\bm{0.83}$(0.22) & 16  \\
& & Joe & 1.09(0.22) & 0.93(0.22) & 1.08(0.22) & 6.38(\underline{0.07}) & 3.43(0.10) & $\bm{0.89}$(0.22) & 20  \\
& & Clayton & 1.13(0.24) & 0.96(0.24) & 1.12(0.25) & 6.49(\underline{0.08}) & 3.49(0.10) & $\bm{0.92}$(0.24) & 20 \\
\hhline{~=========}
& \multicolumn{2}{|c||}{Independent } & 1.77(0.41) & 1.58(0.40) & 1.76(0.41) & 0.68(0.29) & 0.94(0.34) & $\bm{1.17\text{e-}14}$(\underline{0.00}) & 0 \\
\hline
\end{tabular}
}
\caption{\label{MK-SE-copula} Bivariate copulas: Relative MK-SE $\times 100$.
Bold values show the minimun of relative MK-SE and the underline
the minimun of standard deviation for the corresponding copula.}
\end{table}

\FloatBarrier

\begin{table}[!htbp]
\centering
\resizebox{18cm}{8cm}{%
\begin{tabular}{|c|c||c||c||c||c||c||c|}
\hline
\multirow{3}{*}{Copulas} & \multirow{3}{*}{Methods} &  \multicolumn{6}{|c|}{sample size $n=500$}\\\cline{3-8}
& & \multicolumn{3}{c}{$\tau=0.3$}&\multicolumn{3}{|c}{$\tau=0.8$}
\\\hhline{~~======}
 &  & MK-SE $\times 10^3$ & MIAE $\times 10^3$ & MISE $\times 10^5$ & MK-SE $\times 10^3$ &  MIAE $\times 10^5$ &  MISE $\times 10^5$ \\\hhline{========}
\multirow{4}{*}{\rotatebox{90}{Clayton}}
 & Emp & 30.4(5.86) & 5.61(1.39) & 6.28(3.21) & 19.8(4.16) & 1.60(0.27) & 0.85(\underline{0.43})\\[0.5ex]
 & Beta & $\bm{ 27.1}$(\underline{5.42}) & $\bm{5.23}$(1.48) & $\bm{5.57}$(3.26) & $\bm{17.0}$(4.23) & $\bm{1.20}$(0.35) & $\bm{0.70}$(0.49)\\[0.5ex]
 & Check  & 30.4(5.90) & 5.58(\underline{1.38}) & 6.22(\underline{3.20}) & 19.6(4.16) & 1.29(0.27) & 0.80(0.44) \\[0.5ex]
 & CN  & 29.7(7.65) & 6.73(1.40) & 8.35(3.94) & 60.7(\underline{0.25}) & 12.0(\underline{0.06}) & 25.3(0.63)\\[1ex]\hline\hline
 \multirow{4}{*}{\rotatebox{90}{Joe}}
 & Emp & 31.9(6.27) & 5.86(\underline{1.3}) & 6.76(3.26) & 19.8(4.10) & 2.24(0.31) & 1.33(0.46)\\[0.5ex]
 & Beta & 28.3(\underline{5.83}) & 5.50(1.34) & 6.02(\underline{3.11}) & 17.0(3.94) & $\bm{1.83}$(0.37) & $\bm{1.05}$(0.49)\\[0.5ex]
 & Check  & 31.7(6.20) & 5.85(1.31) & 6.74(3.26) & 19.6(4.06) & 2.07(\underline{0.31}) & 4.88(\underline{0.46}) \\[0.5ex]
 & CN  & $\bm{25.7}$(7.20) & $\bm{5.20}$(1.53) & $\bm{5.57}$(3.80) & $\bm{15.5}$(4.12) & 2.10(0.
 460) & 1.10(0.62)\\[1ex]\hline\hline
\multirow{4}{*}{\rotatebox{90}{Gumbel}}
 & Emp  & 30.5(6.28) & 5.80(\underline{1.35}) & 6.48(3.31) & 18.8(3.74) & 2.00(0.27) & 1.09(0.36)\\[0.5ex]
 & Beta & 27.1(5.76) & 5.44(1.43) & 5.76(3.32) & 15.8(3.45) & $\bm{1.60}$(0.33) & 0.83(0.39)\\[0.5ex]
 & Check & 30.2(6.18) & 5.79(1.36) & 6.45(\underline{3.14}) & 18.6(3.76) & 1.82(\underline{0
 .27}) & 1.04(\underline{0.36}) \\[0.5ex]
 & CN  & $\bm{24.3}$(\underline{5.37}) & $\bm{5.06}$(1.48) & $\bm{5.07}$(3.66) & $\bm{14.3}$(\underline{3.49}) & 1.85(0.33) & $\bm{0.82}$(0.37)\\[1ex]\hline\hline
 \multirow{4}{*}{\rotatebox{90}{Frank}}
 & Emp & 29.1(5.15) & 5.5(0.98) & 5.81(2.15) & 17.8(3.32) & 1.82(0.22) & $\bm{1.09}$(0.360)\\[0.5ex]
 & Beta & 25.9(4.68) & 5.16(1.03) & 5.07(2.09) & $\bm{15.5}$(3.15) & 1.42(0.26) & 0.83(0.39)\\[0.5ex]
 & Check & 28.9(5.07) & 5.55(\underline{0.98}) & 5.77(2.14) & 17.7(3.31) & 1.61(0.21) & 1.04(0.36) \\[0.5ex]
 & CN  & $\bm{15.8}$(\underline{3.34}) & $\bm{4.52}$(1.12) & $\bm{3.70}$(\underline{2.04}) & 18.3(\underline{2.12}) & $\bm{0.34}$(\underline{0.17}) & 2.13(\underline{0.24})\\[1ex]\hline\hline
\multirow{4}{*}{\rotatebox{90}{Student}}
 & Emp  & 30.0(6.27) & 5.51(1.034) & 5.81(2.40) & 18.7(\underline{3.97}) & 1.84(0.27) & 0.99(0.36)\\[0.5ex]
 & Beta & 26.9(5.71) & 5.10(1.08) & 5.06(2.35) & 15.6(\underline{3.97}) & $\bm{1.40}$(0.31) & 0.75(0.36)\\[0.5ex]
 & Check & 29.9(6.19) & 5.48(\underline{1.033}) & 5.77(2.40) & 18.7(\underline{3.97}) & 1.61(\underline{0.26}) & 0.95(\underline{0.35}) \\[0.5ex]
 & CN & $\bm{18.6}$(\underline{4.89}) & $\bm{4.40}$(1.25) & $\bm{3.76}$(\underline{2.24}) & $\bm{12.4}$(3.45) & 1.63(0.34) & $\bm{0.64}$(0.38)\\[1ex]\hline\hline
 \multirow{4}{*}{\rotatebox{90}{Normal}}
 & Emp & 30.0(6.01) & 5.46(1.06) & 5.78(2.45) & 18.6(4.16) & 1.81(0.26) & 0.99(0.37)\\[0.5ex]
 & Beta & 26.8(5.80) & 5.05(1.11) & 5.03(\underline{2.42}) & 15.6(4.01) & $\bm{1.34}$(0.32) & 0.74(0.39)\\[0.5ex]
 & Check  & 29.9(6.00) & 5.44(\underline{1.07}) & 5.75(2.47) & 18.5(4.13) & 1.55(\underline{0.26}) & 0.94(\underline{0.36}) \\[0.5ex]
 & CN & $\bm{17.9}$(\underline{5.03}) & $\bm{4.44}$(1.52) & $\bm{3.31}$(3.31) & $\bm{13.8}$(\underline{3.65}) & 1.66(0.33) & $\bm{0.68}$(0.41)\\[1ex]\hline\hline
%
%
\multirow{5}{*}{\rotatebox{90}{Independent}}& &\multicolumn{2}{|c}{MK-SE $\times 10^3$} &\multicolumn{2}{|c}{MIAE $\times 10^3$}& \multicolumn{2}{|c|}{MISE $\times 10^5$}\\\cline{3-8}\cline{3-8}
 & Emp &\multicolumn{2}{|c|}{31.3(5.51)} &\multicolumn{2}{c|}{5.59(1.03)}& \multicolumn{2}{c|}{6.09(2.43)}\\[0.5ex]\cline{3-8}
 & Beta &\multicolumn{2}{|c|}{27.9(5.48)} &\multicolumn{2}{c|}{5.20(1.06)}& \multicolumn{2}{c|}{5.33(2.31)}\\[0.5ex]\cline{3-8}
 & Check &\multicolumn{2}{|c|}{31.1(5.42)} &\multicolumn{2}{c|}{5.57(1.03)}& \multicolumn{2}{c|}{6.05(2.43)}\\[0.5ex]\cline{3-8}
 & CN &\multicolumn{2}{|c|}{$\bm{2.22\text{e-}13}$(\underline{0.00})} &\multicolumn{2}{c|}{$\bm{7.34\text{e-}15}$(\underline{0.00})}& \multicolumn{2}{c|}{$\bm{3.67\text{e-}29}$(\underline{0.00})}\\[1ex]\hhline{========}
\multicolumn{8}{c}{}  \\[0.5ex] \hhline{========}
\multicolumn{8}{c}{parameter optimal N of our procedure}  \\[0.5ex] \hhline{========}
dependence level & Clayton & Joe & Gumbel & Frank & Student ($df=17$) & Gauss & Independent \\\hhline{========}
$\tau=0.3$ & 2 & 9 & 8 & 1 & 3 & 3 & \multirow{2}{*}{0}\\\cline{1-7}
$\tau=0.8$ &          2 & 12 & 12 & 4 & 10 & 13 & \\\hline
\end{tabular}
}
\caption{ \label{trivariate-copula500} Trivariate copulas: MIAE$\times 10^3$,
MISE$\times 10^5$ and MK-SE$\times 10^3$ witth sample size $n=500$. Bold values show the minimun of MIAE, MISE or
MK-SE (it depends on the column) and the underline the minimun of
standard deviation for the corresponding copula.}
\end{table}

\FloatBarrier

\begin{table}[!htbp]
\centering
\resizebox{18cm}{8cm}{%
\begin{tabular}{|c|c||c||c||c||c||c||c||}
\hline
\multirow{3}{*}{Copulas} & \multirow{3}{*}{Methods} & \multicolumn{6}{c|}{sample size $n=200$} \\\cline{3-8}
& & \multicolumn{3}{c}{$\tau=0.3$}&\multicolumn{3}{|c}{$\tau=0.8$} \\\hhline{~~======}
 &  & MK-SE $\times 10^3$ & MIAE $\times 10^3$ & MISE $\times 10^5$ & MK-SE $\times 10^3$ &  MIAE $\times 10^3$ &  MISE $\times 10^5$\\\hhline{========}
\multirow{4}{*}{\rotatebox{90}{Clayton}}
 & Emp & 48.4(9.37) & 8.94(2.47) & 16.1(9.47) & 30.5(6.43) & 2.87(0.39) & 2.71(\underline{0.88}) \\[0.5ex]
 & Beta & $\bm{39.9}$(9.38) & $\bm{8.13}$(2.67) & $\bm{13.5}$(\underline{9.28}) & $\bm{23.9}$(7.13) & $\bm{1.91}$(0.62) & $\bm{1.72}$(1.17) \\[0.5ex]
 & Check & 47.9(9.37) & 8.87(2.45) & 15.8(9.30) & 29.6(7.04) & 1.94(0.39) & 1.87(0.89)  \\[0.5ex]
 & CN & 45.9(\underline{8.16}) & 14.3(\underline{1.73}) & 33.8(9.82) & 61.7(\underline{0.43}) & 0.012(\underline{0.12}) & 25.8(0.93) \\[1ex]\hline\hline
 \multirow{5}{*}{\rotatebox{90}{Joe}}
 & Emp & 49.5(9.39) & 9.06(\underline{1.86}) & 16.0(\underline{7.08}) & 31.4(8.61) & 3.81(0.81) & 3.63(1.94) \\[0.5ex]
 & Beta & 40.2(8.84) & 8.10(2.06) & 13.1(7.10) & 24.8(8.32) & $\bm{2.82}$(0.95) & 2.65(\underline{1.20}) \\[0.5ex]
 & Check & 48.8(9.44) & $\bm{9.01}$(1.90) & 15.8(7.19) & 30.7(8.41) & 3.23(\underline{0.68}) & 3.23(1.62)  \\[0.5ex]
 & CN & $\bm{38.0}$(\underline{8.62}) & 8.04(2.10) & $\bm{13.0}$(7.23) & $\bm{20.5}$(\underline{6.99}) & 3.33(0.90) & $\bm{2.45}$(1.70) \\[1ex]\hline\hline
\multirow{4}{*}{\rotatebox{90}{Gumbel}}
 & Emp & 47.2(8.67) & 8.63(\underline{1.38}) & 14.1(4.92) & 30.4(6.23) & 3.42(0.56) & 3.04(1.18)\\[0.5ex]
 & Beta & 38.4(8.30) & 7.66(1.57) & 11.3(\underline{4.89}) & 23.6(7.29) & $\bm{2.40}$(0.69) & 2.06(1.26) \\[0.5ex]
 & Check & 46.7(8.58) & 8.62(1.40) & 14.0(4.90) & 29.8(5.81) & 2.88(\underline{0.45}) & 2.71(9.61) \\[0.5ex]
 & CN & $\bm{34.6}$(\underline{7.67}) & $\bm{7.41}$ (1.66) & $\bm{11.07}$(4.95) & $\bm{18.4}$(\underline{5.64}) & 2.99(0.59) & $\bm{1.85}$(\underline{0.97}) \\[1ex]\hline\hline
 \multirow{4}{*}{\rotatebox{90}{Frank}}
 & Emp & 47.0(8.97) & 8.91(2.35) & 15.3(9.53) & 29.3(6.85) & 3.25(0.55) & 3.04(1.18) \\[0.5ex]
 & Beta & 38.1(8.65) & 7.94(2.54) & 12.5(9.55) & $\bm{24.8}$(6.61) & $\bm{2.24}$(0.65) & $\bm{2.06}$(1.26) \\[0.5ex]
 & Check & 46.3(8.92) & 8.84(\underline{2.32}) & 15.0(9.42) & 28.5(6.50) & 2.54(0.46) & 2.71(\underline{0.96})  \\[0.5ex]
 & CN & $\bm{19.8}$(\underline{7.20}) & $\bm{5.98}$(2.85) & $\bm{7.29}$(\underline{8.49}) & 56.9(\underline{1.92}) & 10.6(\underline{0.19}) & 20.3(0.97) \\[1ex]\hline\hline
\multirow{4}{*}{\rotatebox{90}{Student}}
 & Emp & 48.6(8.79) & 8.76(\underline{1.66}) & 14.7(6.13) & 29.1(6.13) & 3.15(0.53) & 2.60(0.98) \\[0.5ex]
 & Beta & 39.3(8.23) & 7.76(1.84) & 11.7(6.12) & 22.7(7.25) & $\bm{2.04}$(0.63) & 1.70(1.20)  \\[0.5ex]
 & Check & 47.9(8.90) & 8.73(1.68) & 14.5(6.20) & 28.8(6.49) & 2.44(\underline{0.40}) & 2.25(0.91) \\[0.5ex]
 & CN & $\bm{27.8}$(\underline{6.96}) & $\bm{6.78}$(2.12) & $\bm{9.02}$(\underline{6.05}) & $\bm{16.8}$(\underline{5.18}) & 2.74(0.52) & $\bm{1.49}$(\underline{0.87}) \\[1ex]\hline\hline
 \multirow{4}{*}{\rotatebox{90}{Normal}}
 & Emp & 47.4(7.94) & 8.65(1.51) & 14.4(5.31) & 28.6(5.92) & 3.13(0.43) & 2.58(0.95) \\[0.5ex]
 & Beta & 38.4(7.98) & 7.65(\underline{1.71}) &  11.5(5.51) & 22.5(6.59) & $\bm{1.99}$(0.60) & 1.67(1.08) \\[0.5ex]
 & Check & 46.8(8.19) & 7.65(2.54) & 14.2(5.35) & 28.1(5.92) & 2.35(\underline{0.40}) & 2.21(8.80)  \\[0.5ex]
 & CN & $\bm{24.5}$(\underline{6.28}) & $\bm{6.54}$(1.96) & $\bm{8.08}$(\underline{5.19}) & $\bm{16.3}$(\underline{4.90}) & 2.71(0.50) & $\bm{1.43}$(\underline{0.77})  \\[1ex]\hline\hline
%
%
\multirow{5}{*}{\rotatebox{90}{Independent}}& & \multicolumn{2}{c}{MK-SE  $\times 10^3$} &\multicolumn{2}{c}{MIAE  $\times 10^3$}& \multicolumn{2}{|c|}{MISE  $\times 10^5$} \\\cline{3-8}\cline{3-8}
 & Emp &  \multicolumn{2}{c|}{49.7(8.09)} &\multicolumn{2}{c|}{8.77(1.50)}& \multicolumn{2}{c|}{14.9(5.40)} \\[0.5ex]\cline{3-8}
 & Beta & \multicolumn{2}{c|}{41.4(7.23)} &\multicolumn{2}{c|}{7.82(1.63)}& \multicolumn{2}{c|}{12.0(5.24)} \\[0.5ex]\cline{3-8}
 & Check & \multicolumn{2}{c|}{49.1(8.03)} &\multicolumn{2}{c|}{8.73(1.53)}& \multicolumn{2}{c|}{14.7(5.47)} \\[0.5ex]\cline{3-8}
 & CN & \multicolumn{2}{c|}{$\bm{2.22\text{e-}11}$(\underline{0.00})} &\multicolumn{2}{c|}{$\bm{7.99\text{e-}13}$(\underline{0.00})}& \multicolumn{2}{c|}{$\bm{4.38\text{e-}13}$(\underline{0.00})}
  \\[1ex]\hhline{========}

\multicolumn{8}{c}{}  \\[0.5ex] \hhline{========}
\multicolumn{8}{c}{parameter optimal N of our procedure}  \\[0.5ex] \hhline{========}
dependence level & Clayton & Joe & Gumbel & Frank & Student ($df=17$) & Gauss & Independent \\\hhline{========}
$\tau=0.3$ & 1& 7 & 6 & 1 & 3 & 2 & \multirow{2}{*}{0}\\\cline{1-7}
$\tau=0.8$        & 2 & 7 & 7 & 2 & 7 & 7 & \\\hline
\end{tabular}
}
\caption{ \label{trivariate-copula200} Trivariate copulas: MIAE $\times 10^3$,
MISE $\times 10^5$ and MK-SE $\times 10^3$ with sample size $n=200$. Bold values show the minimun of MIAE, MISE or MK-SE (it depends on the column) and the underline the minimun of
standard deviation for the corresponding copula.}
\end{table}

\FloatBarrier

\begin{table}[!h]
\centering
\resizebox{18cm}{8cm}{%
\begin{tabular}{|c|c|c||c|c|c|cc|c|c|c|}
\hline
\multirow{24}{*}{\rotatebox{90}{sample size $n=500$ }}
&\multicolumn{2}{|c}{\multirow{2}{*}{Copula} }&  \multicolumn{8}{|c|}{Method}\\ \cline{4-11}
&\multicolumn{2}{|c|}{} & $\w{c}_{Mr}$ & $\w{c}_{Bk}$ & $\w{c}_{wa}^{} $ &N &$\w{c}^{[{\mathbf N}]}$ & $\w{c}_{Pt}$ & $\w{c}_{Be10}$ & $\w{c}_{Be25}$\\\hhline{~==========} 
&\multirow{7}{*}{\rotatebox{90}{$\tau=0.3$ }}
& Clayton & 48.3(1.4) & 31.3(3.8) & 68.5(\underline{0.67}) & 5 & $\bm{15.0}$(4.3)  & 23.0(5.4) & 38.7(2.1) & 22.9(4.8)\\  
& & Joe & 63.3(\underline{1.2}) & 43.7(3.7) & 61.7(1.7) & 7 & $\bm{9.1}$(5.0) &  35.9(5.7) & 51.7(2.0)  & 34.2(4.9)\\  
& & Gumbel & 42.6(\underline{1.4}) & 26.9(3.1) & 49.5(1.8) & 3 & 25.2(1.7)  & 19.4(5.0) & 32.3(2.1) & $\bm{19.3}$(3.9)\\
& & Frank & 3.0(1.1) & 2.1(0.75) & 30.4(1.9) & 1 & $\bm{1.9}$(\underline{0.23}) & 2.8(1.2) & 1.9(0.78)  & 2.8(1.0)\\
& & Student(df=17) & 17.0(1.6) & $\bm{9.1}$(2.1) & 38.5(1.5) & 2 & 9.6(\underline{0.87}) & 4.5(2.5) & 10.6(1.6) & 6.2(2.6)\\
& & Gauss & 12.4(1.4) & 6.2(1.8) & 35.9(1.3) & 2 & $\bm{6.2}$(\underline{0.82})  & 3.2(2.0) & 7.1(1.5) & 4.2(2.0)\\
\hhline{~==========}
&\multirow{7}{*}{\rotatebox{90}{$\tau=0.55$ }}
& Clayton & 75.6(0.99) &53.5(3.1) &92.1(0.13) & 9 & $\bm{26.6}$(3.0) &  55.6(5.7) & 74.6(\underline{0.90})  & 58.0(2.0)\\  
& & Joe & 78.6(0.97) & 56.8(2.8) & 79.6(\underline{0.78}) & 10 & $\bm{30.0}$(3.4)  & 58.3(4.7) & 76.6(0.84)  & 61.2(2.1)\\  
& & Gumbel & 63.9(1.4) & 39.1(3.9) & 70.1(\underline{1.1}) & 6 & $\bm{24.2}$(3.6) &  39.9(6.7) & 61.4(1.3)  & 43.4(3.1)\\
& & Frank & 4.2(1.1) & 2.9(0.99) &32.4(1.2) & 3 & $\bm{1.6}$(\underline{0.57})  & 2.5(1.2) & 6.3(1.1)  &2.9(1.1)\\
& & Student & 37.6(1.8) & 16.5(3.6) & 58.3(\underline{1.1}) & 5 & $\bm{9.8}$(2.6) &  14.5(5.2) & 35.7(1.8)  & 18.1(3.3)\\
& & Gauss & 31.6(1.8) & 12.3(3.4) & 54.0(\underline{1.2}) & 5 & $\bm{7.0}$(2.4) & 10.5( 5.1) & 30.0 ( 1.7)  & 13.7(3.4)\\
\hhline{~==========}
&\multirow{7}{*}{\rotatebox{90}{$\tau=0.8$ }}
& Clayton & 83.6(1.1) & 70.3(2.1) & 98.1(\underline{0.02}) & 17 & $\bm{57.7}$(1.0)  & 80.7(3.0) & 91.2(0.18)  & 84.0(0.40)\\  
& & Joe & 84.5(0.92) & 70.6(1.8) & 91.5(\underline{0.19}) & 16 & $\bm{61.5}$(0.87) & 81.4(2.98) & 91.4(0.18)  & 84.3(0.38)\\  
& & Gumbel & 71.1(1.8) & 49.5(3.2) & 85.8(\underline{0.36}) & 10 & $\bm{45.2}$(1.1)   & 67.3(4.6) & 83.6(0.36)  & 71.1(0.79)\\
& & Frank & 6.2(1.7) & 5.4(2.18) & 38.3(\underline{0.65}) & 8 & 5.3(3.0)  & $\bm{4.9}$(2.1) & 28.7(0.91) & 12.9(1.5) \\
& & Student & 50.5(2.5) & 22.4(3.9) & 80.9(\underline{0.46}) & 10 & $\bm{15.4}$(2.5)  &  43.3(7.5) & 70.4(0.68)  & 50.5(1.6)\\
& & Gauss & 45.4(2.7) & 17.3(3.7) & 78.3(\underline{0.48}) & 10 & $\bm{11.0}$(2.2) &   38.3(8.8) & 67.3(0.74) & 46.1(1.7)\\
\hhline{~==========}
&\multicolumn{2}{|c||}{Independent }  & 0.20(0.17) & 0.89(0.38) & 29.8(1.9) & 0 & $\bm{0.0}$(\underline{0.0})  & 2.8(1.2) & 0.89(0.43)  & 3.4(1.0)\\
\hhline{===========}
\multirow{24}{*}{\rotatebox{90}{sample size $n=1000$ }}
&\multicolumn{2}{|c}{\multirow{2}{*}{Copula} }&  \multicolumn{8}{|c|}{Methods}\\ \cline{4-11}
&\multicolumn{2}{|c|}{} & $\w{c}_{Mr}$ & $\w{c}_{Bk}$ & $\w{c}_{wa}$ &N &$\w{c}^{[{\mathbf N}]}$ &  $\w{c}_{Pt}$ & $\w{c}_{Be10}$ & $\w{c}_{Be25}$\\\hhline{~==========}
&\multirow{7}{*}{\rotatebox{90}{$\tau=0.3$ }}
& Clayton & 47.0(1.00) & 27.4(3.20) & 67.1(\underline{0.38}) & 12 & $\bm{6.40}$(2.70)  & 19.8(4.60) & 38.5(1.50)  & 21.8(3.60)\\  
& & Joe & 62.0(\underline{1.00}) & 39.0(3.10) & 60.6(1.20) & 7 & $\bm{5.70}$(2.90)  & 31.4(5.40) & 51.5(1.50) & 32.9(3.10)\\  
& & Gumbel & 41.3(\underline{0.99}) & 23.1(2.40) & 47.1(1.20) & 5 & $\bm{11.2}$(2.40)  & 16.5(4.00) & 31.9(1.40)  & 17.7(3.00)\\
& & Frank & 2.1(0.61) & 1.30(0.45) & 27.2(0.86) & 2 & $\bm{0.46}$(\underline{0.29}) & 1.90(0.66) & 1.40(0.54)  &1.50(0.54)\\
& & Student & 15.7(1.20) & 7.20(1.60) & 35.8(0.88) & 2 & 9.20(\underline{0.54})  & $\bm{3.50}$(1.90) & 10.2(1.30)  & 4.70(1.70)\\
& & Gauss & 11.3(0.88) & 5.00(1.20) & 33.3(0.86) & 3 & 3.80(\underline{0.82})  & $\bm{2.30}$(1.40) & 6.90(1.00)  & 3.30(1.30)\\
\hhline{~==========}
&\multirow{7}{*}{\rotatebox{90}{$\tau=0.55$ }}
& Clayton & 73.8(0.78) & 49.1(2.5) & 91.9(\underline{0.07}) & 13 & $\bm{23.3}$(2.10)  & 50.5(5.70) & 74.5(0.63)  & 56.0(1.6)\\  
& & Joe & 77.1(0.89) & 52.3(2.20) & 79.3(\underline{0.52}) & 17 & $\bm{21.8}$(4.20)  & 53.5(5.90) & 76.5(0.5603)  & 60.7(1.20)\\  
& & Gumbel & 61.8(1.00) & 33.8(2.80) & 69.5(\underline{0.73}) & 12 & $\bm{10.3}$(2.7)  & 35.5(5.80) & 61.1(0.93)  & 42.0(1.90)\\
& & Frank & 3.30(0.68) & 2.00(0.63) & 30.4(0.62) & 3 & $\bm{1.28}$(\underline{0.30})  & 1.90(0.87) & 6.00(0.82) & 2.30(0.71) \\
& & Student(df=17) & 34.9(1.60) & 12.9(2.80) & 56.9(\underline{0.76}) & 5 & $\bm{8.10}$(1.60)  & 12.5(4.50) & 35.5(1.40) & 17.6(2.30) \\
& & Gauss & 29.20(1.40) & 9.70(2.50) & 52.7(\underline{0.80}) & 6 & $\bm{3.50}$(1.40)  & 8.10(3.70) & 29.9(1.20)  & 13.7(2.30)\\
\hhline{~==========}
&\multirow{7}{*}{\rotatebox{90}{$\tau=0.8$ }}
& Clayton & 81.9(1.00) & 67.1(1.80) & 98.1(\underline{0.02}) & 20 & $\bm{46.6}$(1.40)  & 78.4(2.80)  & 9.2(0.12)  & 84.0(0.28)\\  
& & Joe & 83.0(1.00) & 67.9(1.80) & 91.4(\underline{0.14}) & 20 & $\bm{47.5}$(1.2) & 79.1(3.00) & 91.4(0.12)  & 84.3(0.31)\\  
& & Gumbel & 68.2(1.50) & 45.7(2.40) & 85.6(\underline{0.23}) & 16 & $\bm{38.2}$(1.60)  & 63.5(4.90) & 83.5(0.23)  & 70.8(0.53)\\
& & Frank & 4.80(1.10) & 3.90(1.80) & 37.7(\underline{0.42}) & 9 & $\bm{3.6}$(1.7)  &  3.97(1.60) & 28.7(0.65)  & 12.9(1.10)\\
& & Student & 45.8(2.30) & 18.6(3.00) & 80.5(\underline{0.30}) & 13 & $\bm{12.6}$(1.40) & 39.6(7.10) & 70.2(0.46) & 50.4(1.10)\\
& & Gauss & 41.1(2.20) & 14.4(2.70) & 78.1(\underline{0.34}) & 14 & $\bm{8.90}$(1.60) & 32.5(6.70) & 67.2(0.54)  & 46.2(1.20)\\
\hhline{~==========}
&\multicolumn{2}{|c||}{Independent }  & 0.15(0.10) & 0.62(0.22) & 25.90(0.88) & 0 & $\bm{0.00}$(\underline{0.00}) &   1.90(0.64) & 0.44(0.19)  & 1.70(0.51)\\
\hline
\end{tabular}
}
\caption{\label{MISE-density}Bivariate copulas densities: Relative
MISE $\times 100$. Bold values show the minimun of relative MISE and the
underline the minimun of standard deviation for the corresponding
copula.}
\end{table}

\FloatBarrier

\begin{table}[!h]
\centering
\resizebox{18cm}{8cm}{%
\begin{tabular}{|c|c|c||c|c|c|cc|c|c|c|}
\hline
\multirow{24}{*}{\rotatebox{90}{sample size $n=500$ }}
&\multicolumn{2}{|c}{\multirow{2}{*}{Copula} }&  \multicolumn{8}{|c|}{Method}\\ \cline{4-11}
&\multicolumn{2}{|c|}{} & $\w{c}_{Mr}$ & $\w{c}_{Bk}$ & $\w{c}_{wa}$ &N &$\w{c}^{[{\mathbf N}]}$ & $\w{c}_{Pt}$ & $\w{c}_{Be10}$ & $\w{c}_{Be25}$\\\hhline{~==========} 
&\multirow{7}{*}{\rotatebox{90}{$\tau=0.3$ }}
& Clayton & 88.9(1.10) & 71.5(4.30) & 100(\underline{0.00}) & 5 & $\bm{43.4}$(9.20)  & 60.6(7.60) & 79.8(2.20)  & 60.2(6.80) \\  
& & Joe & 92.5(\underline{0.82}) & 76.8(3.30) & 85.6(1.60) & 7 & $\bm{24.8}$(13.3)  & 69.4(5.60) & 83.6(1.70)  & 67.2(4.90) \\  
& & Gumbel & 89.5(\underline{1.10}) & 71.3(4.10) & 81.4(2.30) & 3 & 67.1(3.50)  & 59.9(8.00) & 78.5(2.60)  & $\bm{59.2}$(6.50) \\
& & Frank & 46.5(5.20) & 29.7(9.60) & 100(\underline{0.00}) & 1 & $\bm{23.4}$ (3.30) & 43.4(20.0) & 29.2(10.9) & 44.3(18.3) \\
& & Student & 77.5(2.40) & 58.6(6.90) & 100(\underline{0.00}) & 2 & 60.3(4.20)  & $\bm{38.0}$(12.1) & 64.4(4.60)  & 48.6(14.6) \\
& & Gauss & 73.2(2.70) & 52.6(8.50) & 100(\underline{0.00}) & 2 & 54.4(5.2)  & $\bm{30.5}$(14.7) & 59.0(6.00) & 42.1(15.4) \\
\hhline{~==========}
&\multirow{7}{*}{\rotatebox{90}{$\tau=0.55$ }}
& Clayton & 90.7(0.61) & 76.3(2.10) & 100(\underline{0.00}) & 9 & $\bm{53.6}$(3.10) &   78.0(4.00) & 89.9(0.57)  & 79.4(1.40) \\  
& & Joe & 91.9(0.58) & 78.2(1.90) & 91.4(\underline{0.5}) & 10 & $\bm{56.6}$(3.20)  & 79.4(3.20) & 90.5(0.52)  & 81.1(1.30) \\  
& & Gumbel & 88.3(0.96) & 69.3(3.60) & 87.7(\underline{0.88}) & 6 &$\bm{ 53.2}$(4.40)  & 70.5(5.90) & 86.2(0.97)  & 72.9(2.70) \\
& & Frank & 46.3(5.60) & 29.7(11.4) & 100(\underline{0.00}) & 3 & $\bm{14.4}$(6.40) &   26.1(12.0) & 39.5(6.80)  & 30.6(13.0) \\
& & Student & 77.9(2.00) & 53.4(6.60) & 100(\underline{0.00}) & 5 & $\bm{39.0}$(8.90)  & 49.5(8.40) & 74.8(2.10)  & 56.2(5.90) \\
& & Gauss & 75.3(2.40) & 48.6(7.50) & 100(\underline{0.00}) & 5 & $\bm{33.8}$(10.3)  & 44.0(11.4) & 72.2(2.30)  & 52.0(7.40) \\
\hhline{~==========}
&\multirow{7}{*}{\rotatebox{90}{$\tau=0.8$ }}
& Clayton & 92.6(0.58) & 85.0(1.20) & 100(0.00) & 17 & $\bm{77.1}$(0.69) & 91.0(1.60)  & 96.2(0.10)  & 92.5(\underline{0.22}) \\  
& & Joe & 93.1(0.47) & 85.2(1.10) & 96.5(\underline{0.09}) & 16 & $\bm{79.6}$(0.58) &    91.3(1.60) & 96.3(\underline{0.09})  & 92.6(0.21) \\  
& & Gumbel & 87.8(1.00) & 74.2(2.40) & 94.1(\underline{0.21}) & 10 & $\bm{71.3}$(0.87)  &  85.6(2.60) & 93.6(0.21)  & 87.3(0.53) \\
& & Frank & 0.39.6(6.90) & 31.9(13.6) & 100(\underline{0.00}) & 8 & 32.2(15.1) &  $\bm{23.0}$(9.70) & 56.7(2.40)  & 37.2(7.30) \\
& & Student & 75.5(2.00) & 52.7(4.60) & 100(\underline{0.00}) & 10 & $\bm{43.8}$(4.40)  &  70.2(6.00) & 86.8(0.51)  & 74.5(1.40) \\
& & Gauss & 72.0(2.20) & 47.3(5.30) & 100(\underline{0.00}) & 10 & $\bm{38.2}$(5.10) & 66.5(7.50) & 85.3(0.58)  & 71.7(1.80) \\
\hhline{~==========}
&\multicolumn{2}{|c||}{Independent }  & 17.0(5.80) & 38.5(11.6) & 1.00.7(2.70) & 0 & $\bm{0.00}$(\underline{0.00})  & 133(43.9) & 47.3(16.9)  & 102.7(43.5) \\
\hhline{===========}
\multirow{24}{*}{\rotatebox{90}{sample size $n=1000$ }}
& \multicolumn{2}{|c}{\multirow{2}{*}{Copula} }&  \multicolumn{8}{|c|}{Methods}\\ \cline{4-11}
& \multicolumn{2}{|c|}{} & $\w{c}_{Mr}$ & $\w{c}_{Bk}$ & $\w{c}_{wa}$ &N &$\w{c}^{[{\mathbf N}]}$ &  $\w{c}_{Pt}$ & $\w{c}_{Be10}$ & $\w{c}_{Be25}$ \\\hhline{~==========}
& \multirow{7}{*}{\rotatebox{90}{$\tau=0.3$ }}
& Clayton & 87.9(0.84) & 67.1(3.90) & 100(\underline{000}) & 12 & $\bm{14.4}$(8.80)  & 56.4(7.00) & 79.7(1.50)  & 59.3(5.10)\\  
& & Joe & 91.7(\underline{0.69}) & 72.6(2.90) & 85.4(1.00) & 7 & $\bm{19.7}$(10.5)  & 64.5(5.80) & 83.5(1.20)  & 66.4(3.20)\\  
& & Gumbel & 88.5(\underline{0.85}) & 66.3(3.60)& 80.4(1.50) & 5 & $\bm{40.8}$(7.3)  & 55.8(7.0) & 78.2(1.70)  & 57.6(5.20)\\
& & Frank & 43.8(3.67) & 26.3(8.10) & 100(\underline{000}) & 2 & $\bm{12.4}$(5.40) & 37.0(16.5) & 27.7(9.10)  & 33.5(11.6)\\
& & Student & 75.6(2.10) & 52.5(6.50) & 100(\underline{00}) & 2 & 59.0(2.6)  & $\bm{34.7}$(11.4) & 62.6(4.10)  & 44.2(9.70)\\
& & Gauss & 71.4(2.10) & 48.4(7.40) & 100(\underline{000}) & 3 & 42.1(7.90)  & $\bm{26.9}$(12.8) & 58.2(4.40)  & 40.8(11.7)\\
\hhline{~==========}
& \multirow{7}{*}{\rotatebox{90}{$\tau=0.55$ }}
 & Clayton & 89.6(0.49) & 73.2(1.80) & 100(000) & 13 & $\bm{50.3}$(2.30)  & 74.4(4.20) & 89.8(\underline{0.39})  & 79.4(1.10)\\  
& & Joe & 91.0(0.54) & 75.2(1.50) & 91.3(0.34) & 17 & $\bm{47.4}$(4.90)  & 76.2(4.20) & 90.5(\underline{0.34})  & 80.8(0.79)\\  
& & Gumbel & 86.9(0.69) & 64.6(2.70) & 87.6(\underline{0.63}) & 12 & $\bm{33.3}$(5.00)  & 66.7(5.30) & 86.0(0.67)  & 71.9(1.70)\\
& & Frank & 43.9(4.20) & 26.7(10.2) & 100(\underline{000}) & 3 & $\bm{11.1}$(4.50)  & 25.2(11.4) & 39.1(4.90)  & 29.5(10.5)\\
& & Student & 75.3(1.96) & 47.5(5.80) & 100(\underline{000}) & 5 & $\bm{34.2}$ (6.50)  & 46.3(8.80) & 74.3(1.70)  & 54.9(4.2)\\
& & Gauss & 72.7(2.00) & 43.2(6.60) & 100(\underline{0.00}) & 6 & $\bm{21.2}$(9.10)  & 39.0(9.50) & 71.7(1.70)  & 50.8(5.00)\\
\hhline{~==========}
& \multirow{7}{*}{\rotatebox{90}{$\tau=0.8$ }}
& Clayton & 91.7(0.53) & 83.1(1.10) & 100(\underline{0.00}) & 20 & $\bm{69.4}$(1.10)  & 89.7(1.60) & 96.2(0.07) & 92.5(0.16)\\  
& & Joe & 92.3(0.54) & 83.6(1.10) & 96.5(0.07) & 20 & $\bm{70.1}$(0.96) & 90.1(1.70) & 96.3(\underline{0.06})  & 92.6(0.17)\\  
& & Gumbel & 86.3(0.83) & 71.4(1.90) & 94.1(0.14) & 16 & $\bm{65.4}$(1.40) & 83.4(2.96) & 93.5(\underline{0.13}) & 87.1(0.31)\\
& & Frank & 36.2(6.00) & 27.5(12.4) & 100(\underline{0.00}) & 9 & $\bm{27.4}$(12.9) &  23.4(9.90) & 56.5(1.70)  & 37.5(5.20) \\
& & Student & 72.2(1.70) & 48.0(4.30) & 100(\underline{0.00}) & 13 & $\bm{39.5}$(3.20)  &  67.1(6.20) & 86.6(0.40)  & 74.2(0.97)\\
& & Gauss & 69.0(1.70) & 43.5(4.70) & 78.1(\underline{0.34}) & 14 & $\bm{33.4}$(4.20)  &  61.6(6.20) & 85.2(0.44) & 71.5(1.10)\\
\hhline{~==========}
& \multicolumn{2}{|c||}{Independent }  & 18.1(4.90) & 34.4(11.3) & 100(0.00) & 0 & $\bm{0.00}$(\underline{0.00})  & 115.8(33.5) & 34.2(11.8) & 72.6(25.7) \\
\hline
\end{tabular}
}
\caption{\label{MK-SE-density} Bivariate copula densities :
Relative MK-SE $\times 100$. Bold values show the minimun of relative MK-SE
and the underline the minimun of standard deviation for the
corresponding copula.}
\end{table}

\FloatBarrier

\FloatBarrier

\FloatBarrier

\FloatBarrier

\begin{table}[!h]
\centering
\begin{tabular}{|cc|c|c|c|c|c|c|c|c|c|}
\hline
\multicolumn{2}{|c|}{log returns} & Min & Max & Mean & Median & std.Dev.  & Skewness & Kurtosis &JB P-value  \\\hline
& EUR & -4.204e-02 & 4.735e-02 & -2.000e-07 & -7.351e-05  & 0.006754196 & 0.08255672 & 3.438804 &0 \\  
& GBP & -4.115e-02 & 8.321e-02 & 5.298e-05 & 6.305e-05 & 0.006614447 & 0.76795970 & 9.838880 &0 \\   
& YEN & -4.309e-02  & 5.017e-02 & 2.174e-05 & 0.000e+00 & 0.006867485 & 0.19770234 & 3.835367 &0  \\
& CAD & -6.150e-02 & 6.150e-02 & 2.965e-05 & -7.506e-05 & 0.006385885 & 0.05111347 & 13.485304 &0  \\
& AUD & -7.808e-02  & 1.015e-01 & -3.152e-05 & -1.840e-04 & 0.008720500 & 0.58507427  & 11.130791 &0  \\
& CHF & -8.390e-02 & 1.547e-01 & 8.765e-05 & 0.000e+00  & 0.007625693 & 1.64779228 & 43.205732 &0 \\\hline
\end{tabular}
\caption{\label{sum_data} Descriptive statistics for daily log returns}
\end{table}

\FloatBarrier

\begin{table}[!h]
\centering
\begin{tabular}{|c|c|c|c|c|c||c|c|c|}
\hline
cop &$n$ &$\tau$ &  \multicolumn{2}{|c|}{MK-SE} & \multicolumn{2}{|c|}{MISE}  & \multicolumn{2}{|c|}{$N$} \\\hhline{~~~======}\cline{1-3}
\multirow{6}{*}{\rotatebox{90}{Clayton }}&\multirow{3}{*}{\rotatebox{90}{ $500$ }}
& $0.3$ & \multicolumn{2}{|c|}{16.47(10.57)} &  \multicolumn{2}{|c|}{3.98( 4.51)} & \multicolumn{2}{|c|}{4}  \\  
& & $0.55$ & \multicolumn{2}{|c|}{59.39(5.53)} & \multicolumn{2}{|c|}{34.44(6.56)} & \multicolumn{2}{|c|}{10} \\  
& & $0.8$ & \multicolumn{2}{|c|}{61.10(2.90)} & \multicolumn{2}{|c|}{36.54(3.54)} & \multicolumn{2}{|c|}{18} \\
\hhline{~========}
&\multirow{3}{*}{\rotatebox{90}{ $1000$ }}
& $0.3$ & \multicolumn{2}{|c|}{10.12(5.99)} & \multicolumn{2}{|c|}{1.60(1.69)} & \multicolumn{2}{|c|}{3} \\  
& & $0.55$ & \multicolumn{2}{|c|}{60.47(3.73)} & \multicolumn{2}{|c|}{35.40( 4.40)} & \multicolumn{2}{|c|}{11} \\  
& & $0.8$ & \multicolumn{2}{|c|}{57.25(2.46)} & \multicolumn{2}{|c|}{31.98(2.79)} & \multicolumn{2}{|c|}{20}  \\
\hhline{=========}
%
%
\multirow{6}{*}{\rotatebox{90}{Gumbel }}&\multirow{3}{*}{\rotatebox{90}{ $500$ }}
%
%
&  $0.3$ & \multicolumn{2}{|c|}{67.59( 4.22)} & \multicolumn{2}{|c|}{ 39.77(4.91)} & \multicolumn{2}{|c|}{6} \\  
& & $0.55$ & \multicolumn{2}{|c|}{52.70(4.52)} & \multicolumn{2}{|c|}{ 25.52(4.36)} & \multicolumn{2}{|c|}{3} \\ 
& & $0.8$ & \multicolumn{2}{|c|}{ 71.33(0.92)} & \multicolumn{2}{|c|}{ 46.31(1.24)} & \multicolumn{2}{|c|}{10} \\ 
\hhline{~========}
&\multirow{3}{*}{\rotatebox{90}{ $1000$ }}
& $0.3$ & \multicolumn{2}{|c|}{54.36( 5.29)} & \multicolumn{2}{|c|}{ 25.91(4.95)} & \multicolumn{2}{|c|}{4} \\  
& & $0.55$ & \multicolumn{2}{|c|}{33.02(5.58)} & \multicolumn{2}{|c|}{10.64( 3.53)} & \multicolumn{2}{|c|}{11} \\  
& & $0.8$ & \multicolumn{2}{|c|}{68.81(0.80)}  & \multicolumn{2}{|c|}{42.78(1.04)} & \multicolumn{2}{|c|}{14} \\ 
\hline
\end{tabular}
%
\caption{\label{modifErros} Bivariate copula
densities using shrinkage function: Relative MK-SE and
relative MISE in percentages and in brackets the
minimun of standard deviation for the corresponding
copula.
}
\end{table}

\FloatBarrier

\begin{table}[!h]
\centering
\begin{tabular}{|c|c|c|c|c|c|}
\hline
& EUR & GBP & YEN & CAD  & AUD \\\cline{2-6}
 GBP & 10 & & & & \\   
YEN & 3 & 5 &  &  &  \\
 CAD & 8 & 7 & 7 &  & \\
AUD & 7 & 10 & 5 & 11 &\\
 CHF & 16 & 9 & 7 & 5 & 5\\
\hline
\end{tabular}
\caption{\label{parameter_data} Optimal parameter $N$ for daily log-returns}
\end{table}

\FloatBarrier
\begin{table}[!h]
\centering
\begin{tabular}{|c|c|c|c|c|c|}
\hline
 & EUR & GBP & YEN & CAD  & AUD  \\\cline{2-6}
 GBP & 0.4367622 & & & & \\   
YEN & -0.148341452 & -0.119612388 &  &  &  \\
 CAD & -0.182688503 & -0.163666875 & 0.006203468 &  &  \\
 AUD & 0.2668624 & 0.2555104 & -0.1018401 & -0.1843126 &  \\
CHF & -0.6308960 & -0.3978123 & 0.2055829 & 0.1277914 & -0.2176760 \\
\hline
\end{tabular}
\caption{\label{kendall_data} Kendall's rank correlations for daily log-returns}
\end{table}

\FloatBarrier

\FloatBarrier

\FloatBarrier

\FloatBarrier

\end{document}